\documentclass[10pt]{article}
\usepackage{tikz}
\usetikzlibrary{decorations.markings}
\usetikzlibrary{arrows}
\usetikzlibrary{arrows.meta}
\usepackage{amssymb}
\usepackage{amsfonts}
\usepackage{amsmath}
\usepackage{amsthm}
\usepackage{graphicx}
\usepackage{epsfig}
\usepackage{psfrag}
\usepackage{color}
\usepackage{dsfont}

\makeatletter
\newtheorem{lemma}{Lemma}[section]
\newtheorem{theorem}[lemma]{Theorem}

\newtheorem{example}[lemma]{Example}
\newtheorem{definition}[lemma]{Definition}

\newtheorem{proposition}[lemma]{Proposition}
\newtheorem{remark}[lemma]{Remark}
\newcommand{\mR}{\mathbb{R}}
\newcommand{\mL}{\mathbb{L}}
\newcommand{\ML}{\mathds{L}}
\newcommand{\mP}{\mathbb{P}}
\newcommand{\mZ}{\mathbb{Z}}
\newcommand{\J}{\mathcal{J}}

\newcommand{\T}{\mathcal{T}}
\newcommand{\X}{\mathcal{X}}
\newcommand{\Z}{\mathcal{Z}}
\newcommand{\W}{\mathcal{W}}
\newcommand{\cL}{\mathcal{L}}

\newcommand{\cT}{\top}

\newcommand{\Exp}{\mathrm{\,Exp}}

\newcommand{\bq}{\begin{equation}}
\newcommand{\eq}{\end{equation}}
\newcommand{\nameref}{\ref}

\def\BibTeX{{\rm B\kern-.05em{\sc i\kern-.025em b}\kern-.08em
    T\kern-.1667em\lower.7ex\hbox{E}\kern-.125emX}}

\title{\LARGE \bf Classical Thermodynamics Revisited: \\
A Systems and Control Perspective}

\author{Arjan van der Schaft
\thanks{A.J. van der Schaft is with the Bernoulli Institute for Mathematics, Computer
Science and AI, and the Jan C. Willems Center for Systems and Control, University of Groningen, PO Box 407, 9700 AK, the
Netherlands
        {\tt\small a.j.van.der.schaft@rug.nl}}
}

\begin{document}

\maketitle
\thispagestyle{empty}
\pagestyle{empty}

\section{Introduction}
Let us start with the following famous quote from Albert Einstein's autobiographical notes \cite{einstein}, where he expressed his admiration for the theory of classical, macroscopic, thermodynamics:

{\it A theory is more impressive the greater the simplicity of its premises, the more different things it relates, and the more extended its area of applicability. Hence the deep impression that classical thermodynamics made upon me.
It is the only physical theory of universal content concerning which I am convinced that, within the framework of the applicability of its basic concepts, it will never be overthrown.}

On the other hand, there are other, and less favorable, opinions. The eminent Russian mathematician Vladimir I. Arnold stated in \cite{arnold-gibbs}:

{\it Every mathematician knows that it is impossible to understand any elementary course in thermodynamics.
The reason is that the thermodynamics is based, - as Gibbs has explicitly proclaimed -, on a rather complicated mathematical theory, on the contact geometry.}

While Arnold attributed the difficulty in understanding thermodynamics to an inherent mathematical difficulty, others often criticized the imprecise and mathematically outdated presentation of thermodynamics. The perhaps most salient opinion in this regard was expressed by the American mathematician and natural philosopher Clifford Truesdell, who described in \cite[p.6]{truesdell-comical} the formulation of the theory of thermodynamics as a 'dismal swamp of obscurity'. Quoting from \cite{willemsreview}: in a desperate attempt to try to make sense of the writings of De Groot, Mazur, Casimir, and Prigogine, he (Truesdell) goes on to state that there is \cite[p.134]{truesdell} 'something rotten in the [thermodynamic] state of the Low Countries'. 
Clearly, the author of the present paper feels addressed by this last statement.

So much for setting the stage. The purpose of this paper is two-fold. First, to make clear (and de-mystify) the basic concepts of classical thermodynamics, and thus to enable the integration of thermodynamics within systems modeling and control. Second, to demonstrate that systems and control theory provides a natural context for the formulation and understanding of classical thermodynamics. This is not so surprising since classical thermodynamics, firmly rooted in engineering with questions such as the maximal efficiency of steam engines, deals from the very start with systems in interaction with their surrounding (by heat flow, mechanical work, flow of matter, etc.). In particular, it will be shown that dissipativity theory (as founded by Willems in \cite{willems1972}, already with thermodynamics as one of its motivations) is key in the formulation and interpretation of the First and Second Law of thermodynamics. Also a {\it geometric} view on the state properties and the dynamics of thermodynamic systems will be emphasized, thereby unifying and simplifying different representations of thermodynamic systems. On the other hand, I will also argue that thermodynamics motivates paradigm shifts within systems and control; in particular, the use of non-minimal state space formulations, and a geometric view on them. Furthermore, while systems and control theory has been primarily based on linear systems with quadratic cost criteria, in line with basic system models in electrical and mechanical engineering ($RLC$-circuits, mass-spring-damper systems, etc.), thermodynamics necessitates to go beyond this linear-quadratic paradigm. 
From an applications point of view, dealing with thermal behavior is becoming more and more important in many advanced engineering problems. Although the present paper does not address such problems, I hope it paves the way for incorporating, in a natural and unified way, thermodynamics in modeling and control frameworks for solving these problems. 

\subsection{What this paper is {\it not} about}
Thermodynamics is a theory relating many (if not all) areas in the physical sciences, from gases, chemical reactions, to radiation \cite{kondepudi, zemansky}, in accordance with the above quote by Albert Einstein on the universality of thermodynamics. The present paper is {\it not} about all these different application areas, and the presented examples are very simple and do {\it not} illustrate the power and intricacies of the theory when applied to complex physical situations. Instead this paper confines itself to the conceptual and mathematical structure of the theory of classical macroscopic thermodynamics, from the point of view of systems and control theory, as well as of geometry.

Furthermore, the paper concentrates on macroscopic thermodynamics, without making any connection to statistical thermodynamics, or middle-ground theories as in \cite{haddad}. In some sense this is the beauty and power of classical, macroscopic, thermodynamics: it can be presented as a self-contained theory purely based on macroscopic quantities, and the postulates of the First and Second Law. Only in the sidebar "\nameref{otherentropy}" I will briefly point out some connections to other definitions of entropy.
The paper is, moreover, completely about the lumped-parameter case, although much is extendable to mesoscopic (that is, in between macroscopic and microscopic) thermodynamics. Also, I have not made an attempt to provide a balanced view on the large amount of literature about the subject; see e.g. \cite{haddad, kondepudi, zemansky} for a wealth of literature. Finally, the paper reflects my personal journey in understanding classical thermodynamics from a systems and control perspective, and as such presents a biased view, which I hope will stimulate others.

\subsection{Basic terminology} 
\label{subsec:basic}
One of the possible obstacles in understanding classical thermodynamics is the terminology that is being used. This subsection is aimed at explaining, in an informal manner and without claiming full correctness, some of this terminology to a systems and control audience.

First of all, a {\it closed} thermodynamic system is a system that does not exchange matter, but {\it can} exchange energy (through heat or work) with its surrounding. If it can exchange both energy and matter it is called an {\it open} system. A closed thermodynamic system that also does not exchange energy is called an {\it isolated} system. Thus a gas confined in a closed vessel whose volume is controlled by a piston is a {\it closed} system, while it is an {\it isolated} system if it is thermally isolated and the piston does not perform work on the surrounding. Note that this use of terminology is different from the one in systems and control theory, where we would call an isolated thermodynamic system a closed system, and an open system if there is exchange with the surrounding either by matter or energy flow. 

Variables that are proportional to the amount of matter are called {\it extensive} variables, and their value for the overall system is the summation of their values for the parts into which the system is divided. Extensive variables thus scale with the size of the system. Examples of extensive variables are volume, energy, and mass. In contrast, {\it intensive} variables are not matter dependent in nature, and do not depend upon on the total size of the system. Examples of intensive variables are temperature, pressure, and chemical potentials. Note that there are close analogies with e.g. electrical network theory: the charge in a capacitor could be called an extensive variable, while the voltage is an intensive variable. Similarly, the momentum of a point mass would be an extensive variable, while velocity is rather an intensive variable. There is also some relation with the use of 'through' and 'across' variables in physical systems modeling: across variables can be called intensive variables, while time-integrals of through variables are extensive variables.

Although intensive variables such as temperature and pressure may spatially vary within the thermodynamic system, we will restrict attention to thermodynamic systems where the intensive variables are either spatially constant, or the spatial domain can be split into a finite number of parts on each of which the intensive variables (in particular the temperature) are constant. The first situation is often referred to as a homogeneous system, and the second as a nonhomogeneous, or composite, system (consisting of a number homogeneous parts). A state of equilibrium is characterized by the temperature and other intensive variables being uniform throughout the system. If the temperature (and other intensive variables) are not uniform but are well defined locally this is also referred to as 'local equilibrium'. Whenever intensive variables are varying with the spatial position, distributed-parameter (partial differential equations) models are required. This is outside the scope of the present paper, dealing with finite-dimensional, lumped-parameter, thermodynamic systems.

Thermodynamic systems whose intensive variables are a single temperature $T$ and a pressure $P$ (thus taking a uniform value at each time) will be referred to in this paper as {\it simple thermodynamic systems}. A typical example is a single-constituent gas in a container with volume $V$ (an extensive variable of the system), where the variables $V,P,T$ are related through an equation $f(V,P,T)=0$, called the {\it equation of state}. For example, an ideal gas satisfies the equation of state $PV=NRT$, with $N$ the number of moles of the gas and $R$ the gas constant.

Thermodynamic processes or transformations of a simple thermodynamic system are the conversion of one state (i.e., a triple $(V,P,T)$ satisfying the equation of state) to another state. Thermodynamic processes are typically due to interaction with the surrounding (e.g., a piston changing the volumes of a gas container, absorption of heat from an external heat source), but may be also due to an internal local equilibrium (e.g., two heat compartments with different temperatures connected by a conducting wall). A main source for confusion is the use of terminology in this context like 'quasi-reversible', 'infinitesimally slow', or more recently, 'horse-carrot' transformations. As discussed in Remark \ref{rem:carnot} in the section "\nameref{sec:secondlaw}" such terminology can be avoided from a systems and control point of view.

\subsection{Notation}
Notation is fairly standard. Given a function $H:\X \to \mR$ for some $n$-dimensional manifold $\X$ (e.g., $\mR^n$) we will denote by $\frac{\partial H}{\partial x}(x)$ the $n$-dimensional {\it column} vector of partial derivatives, and by $\frac{\partial H}{\partial x^\top}(x)$ the $n$-dimensional {\it row} vector of partial derivatives. Likewise, vectors $v$ will be column vectors, with $v^\top$ denoting its transpose (a row vector).

\section{The First Law of Thermodynamics}
\label{sec:firstlaw}
The First Law of thermodynamics basically expresses two fundamental properties: (1) the different types of interaction of a thermodynamic system with its surrounding, e.g., heat flow, mechanical work, flow of chemical species, etc., all result in an exchange of a common quantity called {\it energy}, (2) there exists a function of the state of the thermodynamic system that represents the {\it internal energy} stored in the system, and which is such that the increase of this function during any time-interval is equal to the sum of the energies supplied to the system by the different forms of interaction with the surrounding during this time-interval ({\it conservation of energy}). Thus energy may manifest itself in different physical forms, which are {\it equivalent} and to a certain extent {\it exchangeable}. 'To a certain extent' because a thermodynamic system can not freely convert one form of energy into another. In particular there are {\it restrictions} to this conversion expressed by the Second Law of thermodynamics, as discussed in the section "\nameref{sec:secondlaw}".

Although the concept of energy at this point in the history of science may seem evident, one should not underestimate the leap in abstraction which lies behind the formulation of the First Law. Energy can{\it not} be directly measured, unlike macroscopic thermodynamic quantities such as volume, pressure and temperature. In particular, the formulation of the equivalence of heat with other forms of energy took form rather late; around the mid of the 19th century, slowly replacing the caloric theory of heat fluid. For a brief historical perspective on the birth of the First Law see e.g. Chapter 2 of \cite{kondepudi}, and the sidebar "\nameref{sb:caloric}".

How to express the First Law in a precise mathematical formulation? Consider a simple thermodynamic system, described by volume $V$, pressure $P$ and temperature $T$.
The mechanical power (rate of mechanical work) provided by the surrounding to the thermodynamic system is given by
\bq
-Pu_V = \mbox{ rate of mechanical work},
\eq
where $u_V:=\dot{V}$ is the rate of volume change. (In order to stick to the usual notation in thermodynamics we follow for the pressure $P$ the physics convention, so that $Pu_V$ is the rate of mechanical work exerted by the system {\it on} the surrounding.) Thus the mechanical work done by the surrounding on the system during a time-interval $[t_1,t_2]$ is
\bq
- \int_{t_1}^{t_2}P(t)u_V(t) dt = - \int_{t_1}^{t_2} P(t)dV(t)
\eq
Second type of interaction with the surrounding is by {\it heat} delivered to the system from a {\it heat source}. Let us denote by $q$ the heat flow (heat per second) from the heat source into the system.
Then the {\it First Law} is expressed by the existence of a function $E(x)$ of the thermodynamic state $x$ (e.g., $(V,P,T)$ satisfying the equation of state), such that along all the possible trajectories of the thermodynamic system
\bq
\label{firstlaw}
E(x(t_2)) - E(x(t_1))  =  \int_{t_1}^{t_2}  \big[q(t) - P(t)u_V(t)\big] dt 
\eq
for all initial conditions $x(t_1)$ and all $t_1 \leq t_2$.
That is, the increase of the total energy $E$ of the thermodynamic system is equal to the incoming heat flow (through the thermal port) minus the mechanical work performed by the system on its surrounding (through the mechanical port). 

Cyclo-dissipativity theory (as explained in the sidebar "\nameref{sb:cyclo}"), immediately yields the right formalism to express the First Law of thermodynamics. Namely, the First Law amounts to the system being {\it cyclo-lossless} for the supply rate $s(q,P,u_V)=q - Pu_V$, with storage function $E$. In fact, we could equally start from cyclo-losslessness with respect to $x^*$, and then infer $E$ as the unique (up to a constant) storage function. Furthermore, in case $E$ is bounded from below (and thus can be turned into a non-negative storage function by adding a suitable constant) the thermodynamic system is {\it lossless}.

The formulation of the First Law is directly extended from simple thermodynamic systems to more involved ones. For example, consider the situation that apart from mechanical and thermal interaction with the surrounding, there is additional mass inflow of chemical species. Then the supply rate $q - Pu_V$ needs to be extended to $q - Pu_V + \sum_k \mu_k \nu_k$. Here $\nu_k=\frac{dN_k}{dt}$, with $N_k$ the mole number of the $k$-th chemical species, and $\mu_k$ its chemical potential; see also \cite{openchemical}.

By using energy as the {\it lingua franca} between different physical domains (mechanical, thermal, electrical, chemical, ..), the First Law is at the heart of the modeling of complex multi-physics systems. The First Law also emphasizes the role of multi-physics systems as {\it energy-converting} devices; energy from one physical domain is converted into energy in another domain. 'Maximal' conversion of heat into mechanical work, motivated by the design of steam engines, was one of the starting points of thermodynamic theory. Electro-chemical devices such as batteries, and electro-mechanical systems including electrical motors and generators, are among the many other examples, often dating back to the nineteenth century or even before \cite{kondepudi}. A second wave of research interest in multi-physics systems based on energy exchange was initiated around the mid of the twentieth century, with an emphasis on the development of unified mathematical modeling and simulation languages for multi-physics systems. In particular, this led to the theory of port-based modeling, bond graphs \cite{paynter, golo}, and eventually port-Hamiltonian systems \cite{maschkebordeaux, vanderschaftmaschkearchive, schjeltsema,passivitybook}. 

However, there is more. While the First Law emphasizes the (lossless) conversion of one form of energy into another, it was soon realized, in fact from the very start of the development of thermodynamic theory, that there are intrinsic {\it limitations} to this energy-conversion. In particular, heat can{\it not} be just converted into mechanical work. This is the origin of the {\it Second Law of thermodynamics}. In classical terminology, while the First Law prohibits the existence of a perpetuum mobile of the first kind (energy cannot be created), the Second Law prohibits the existence of a perpetuum mobile of the second kind (heat cannot be freely converted into mechanical work). This is the topic of the next section.

\section{The Second Law of Thermodynamics} 
\label{sec:secondlaw}

The cyclo-dissipativity interpretation of the Second Law of thermodynamics is much more involved than that of the First Law. Let us start with the formulation of the Second Law as given by Lord Kelvin (see \cite{fermi}):

{\it A transformation of a thermodynamic system whose only final result is to transform into work heat extracted from a source which is at the same temperature throughout is impossible. }

Since the work done during a time-interval $[t_1,t_2]$ is equal to $\int_{t_1}^{t_2}-P(t)dV(t) = \int_{t_1}^{t_2}-P(t)u_V(t)dt$, where, as before, $\dot{V}=u_V$ is the rate of volume change, Kelvin's formulation immediately implies that for each constant temperature $T$ any thermodynamic system is {\it cyclo-passive} with respect to the supply rate $-Pu_V$. 
However, the Second Law is {\it stronger} than just cyclo-passivity for each constant $T$. 
Namely, Kelvin's formulation also forbids the conversion into work of heat from a source at constant temperature for all transformations in which the system interacts as well with a {\it second} heat source at {\it another} temperature, as long as the net heat taken from this second heat source is zero. 

\subsection{Carnot cycle}
\label{subsec:carnot}
In fact, the interaction with heat sources at different temperatures is crucial, as demonstrated by the famous {\it Carnot cycle} (due to Sadi Carnot, 1824). It can be described as follows. Consider a simple thermodynamic system, in particular, a fluid or gas in a confined space of a certain volume. Control the system in two ways: (1) via {\it isothermal} transformations, where {\it heat} is supplied to, or taken from, the system at a constant temperature (classically described as the interconnection of the thermodynamic system with an infinite heat reservoir at the temperature of the isothermal process), (2) via {\it adiabatic} transformations, where the only interaction with the surrounding is via {\it work} supplied to, or taken from, the system (classically described by the movement of a piston that changes the volume of the system, with a pressure equal to the pressure of the gas). Note that while during adiabatic transformations there is no heat absorbed or expelled (but only mechanical work), in isothermal transformations the thermodynamic system is interacting with the surrounding both by heat and by mechanical work (but in such a way that the temperature remains constant).

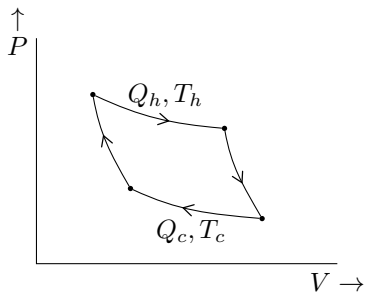
\begin{figure}[h]
\centering
\begin{tikzpicture}

\tikzstyle{bullet} = [circle, fill,minimum size=2pt, inner sep=0pt, outer sep=0pt];

\draw  (1,3) -- (1,0) -- (5,0);
\node at (5,-.25) {$V \rightarrow$};
\node at (0.75,2.90) {$P$};
\node at (0.75,3.25) {$\uparrow$};

\node[bullet] (v4) at (1.75,2.25) {};
\node[bullet] (v3) at (2.25,1) {};
\node[bullet] (v2) at (4,0.6) {};
\node[bullet] (v1) at (3.5,1.8) {};

\begin{scope}[decoration={
    markings,
    mark=at position 0.6 with {\arrow[scale=1.2]{angle 60}}}
    ] 
\draw[-,postaction={decorate}] (v1) to [out=-80,in=120] (v2);
\draw[-,postaction={decorate}] (v2) to [out=175,in=-25] (v3);
\draw[-,postaction={decorate}] (v3) to [out=120,in=-80] (v4);
\draw[-,postaction={decorate}] (v4) to [out=-25,in=175] (v1);
\end{scope}

\node at (2.7151,2.2397) {$Q_h, T_h$};
\node at (3.0541,0.4282) {$Q_c, T_c$};

\end{tikzpicture}
\caption{The Carnot cycle}
\label{fig:carnot}
\end{figure}

A Carnot cycle consists of two isothermal transformations and two adiabatic transformations: first apply an isothermal transformation at temperature $T_h$ ('hot') taking the system from an initial state to another state, secondly apply an adiabatic transformation lowering the temperature of the system to $T_c$ ('cold'), thirdly an isothermal transformation at temperature $T_c$ taking the system to a state from which, fourthly, an adiabatic transformation takes the system back to the original initial state; see Figure \ref{fig:carnot}.

Denote the heat supplied to the system during the first isothermal transformation (at temperature $T_h$) by $Q_h$, and the heat supplied to the system during the other isothermal transformation (at temperature $T_c$) by $Q_c$ (in case heat is actually released from the system, this means $Q_c$ is negative). Such a cycle will be denoted by $(Q_h,T_h; Q_c, T_c)$. For the cycle $(Q_h,T_h; Q_c, T_c)$ being a {\it Carnot cycle} it is moreover assumed that also $(-Q_c,T_c; -Q_h, T_h)$ is a feasible cycle as well. Thus a Carnot cycle can be {\it reversed}, leading to the same initial state. Classically this is referred to as a {\it reversibility}, although this is often understood in a stronger sense: the path in the state space of the thermodynamic system resulting from applying the four transformations can be also followed in the reverse direction; also \cite{zemansky}. In our relaxed definition of 'reversibility', however, we only require that the amounts of heat supplied or released is the negative of those of the original cycle. We will call this {\it heat supply reversibility}. This is summarized as follows.
\begin{definition}
Consider a simple thermodynamic system. A {\it cycle} $(Q_h,T_h; Q_c, T_c)$ from an initial state consists of the alternation of an isothermal (with temperature $T_h$ and supplied heat $Q_h$), an adiabatic, an isothermal (with temperature $T_c$ and supplied heat $Q_c$), and a final adiabatic transformation, which returns the system to the initial state. It is called a Carnot cycle if there also exists a cycle $(-Q_c,T_c; -Q_h, T_h)$.
\end{definition}
\begin{remark}
\label{rem:carnot}
Unfortunately, the exposition of the Carnot cycle is often obscured by the use of terminologies such as 'infinitesimally slow', 'quasi-reversible', 'such that the system remains at equilibrium', etc.. This is largely due to a scientific tradition that thinks about interaction of the system with its surrounding in a different way than in, e.g., systems and control. For example, an isothermal transformation is regarded as the result of the 'real' physical action of a force exerted by a piston on the gas (implying that the pressure delivered by the piston could be different from the pressure of the gas). Furthermore, during the time the system is actuated by the piston it is considered to be in 'real' physical contact with a heat reservoir at a certain temperature (and thus the temperature of the heat reservoir could differ from that of the gas). In contrast, within systems and control (or, e.g., electrical network) theory we have become used to the concept of an 'ideal' control action, where e.g. the pressure and the temperature are directly controlled. This idealized systems point of view is very appropriate for the description of the Carnot cycle. There is a good reason for this; the fundamental problem that was in the mind of Carnot was to give an answer to the {\it control} problem of optimal conversion of heat into work.
\end{remark} 
Importantly, the Second Law can be shown to imply that for Carnot cycles $(Q_h,T_h; Q_c, T_c)$ the quantity
\bq
\frac{Q_h}{-Q_c}
\eq
only depends on the temperatures $T_h,T_c$. This will have some crucial consequences, leading to the celebrated Clausius (in)equality, and eventually to the definition of entropy. In the following we will indicate the main line of reasoning, based on the nice exposition (in a slightly different setting) in \cite{fermi}.

First we start with the following observation. Let us consider a cycle $(Q_h,T_h; Q_c, T_c)$ from a given initial state. Assume that the mechanical work $W := \int P(t)dV(t)$ done {\it by} the thermodynamic system on the surrounding during this cycle is {\it positive}, that is, $W>0$. We claim that this implies $Q_h>0, Q_c<0$ (i.e., the thermodynamic system absorbs heat during the isothermal transformation corresponding to the high temperature $T_h$, and expels heat during the the isothermal transformation corresponding to the low temperature $T_c$). Assume on the contrary that $Q_c\geq 0$. Then we invoke a very reasonable, but in principle {\it additional}, assumption: 

{\it If we bring two heat reservoirs with temperatures $T_h > T_c$ into contact, then any positive amount of heat can be disposed from the reservoir with temperature $T_h$ to that with temperature $T_c$.}

Now complement the cycle $(Q_h,T_h; Q_c, T_c)$ with the flow of an amount of heat from the hot ($T_h$) to the cold ($T_c$) reservoir which is equal to $-Q_c$. Then the net amount of heat absorbed by the cold reservoir is zero, and thus by Kelvin's formulation of the Second Law the mechanical work $W$ satisfies $W\leq 0$, yielding a contradiction. Thus necessarily $Q_c<0$. Furthermore, by the First Law $W=Q_c + Q_h$. Hence, since $W>0$ and $Q_c<0$, it follows that $Q_h>0$ as claimed.

\begin{remark}
Consequently, for a Carnot cycle $(Q_h,T_h; Q_c, T_c)$ with $W>0$, the reversed cycle $(-Q_c,T_c; -Q_h, T_h)$ satisfies $-Q_c>0, -Q_h<0$. Furthermore, the work done by the system on the surrounding during this reversed cycle equals $-Q_c -Q_h =-W<0$. This corresponds to a {\it refrigerator} or {\it heat pump}, where mechanical work is done on the system in order to transfer heat from the cold reservoir to the hot reservoir.
\end{remark}

Now let $(Q_h,T_h; Q_c, T_c)$ be a Carnot cycle, and consider another cycle $(Q'_h,T_h; Q'_c, T_c)$ (from the same initial state); also with $Q_h'>0,Q'_c<0$. First assume that the fraction $\frac{Q_h}{Q'_h}$ is a non-negative rational number, i.e., 
\bq
\label{rat}
\frac{Q_h}{Q'_h} = \frac{N'}{N},
\eq
for some non-negative integers $N,N'$. Consider the total transformation consisting of $N$ cycles $(-Q_c,T_c; -Q_h, T_h)$ and $N'$ cycles $(Q'_h,T_h; Q'_c, T_c)$. Then the total amount of heat absorbed from the hot reservoir $T_h$ is by construction zero. Hence by Kelvin's formulation, the total work $W$ done by the system satisfies $W\leq 0$. On the other hand, by the First Law $W= N(-Q_c) +N'Q'_c$, and thus $N'Q'_c \leq NQ_c$. Together with \eqref{rat} this yields
\bq
\label{carnotin}
\frac{Q_h}{-Q_c} \geq \frac{Q'_h}{-Q'_c}
\eq
Since any real number $\frac{Q_h}{Q'_h}$ can be approximated arbitrarily well by a rational number, this inequality holds for {\it any} Carnot cycle $(Q_h,T_h; Q_c, T_c)$ and any other cycle $(Q'_h,T_h; Q'_c, T_c)$ from the same initial state with $Q_h'>0,Q'_c<0$. 
Furthermore, if $(Q'_h,T_h; Q'_c, T_c)$ is a Carnot cycle as well, then by exchanging the two in the above reasoning, we obtain the opposite inequality, thus proving
\bq
\label{carnote}
\frac{Q_h}{-Q_c} = \frac{Q'_h}{-Q'_c}
\eq
Hence the fraction $\frac{Q_h}{-Q_c}$ is the same for all Carnot cycles between the temperatures $T_c,T_h$, and thus
\bq
\frac{Q_h}{-Q_c}=f(T_c,T_h)
\eq
for some function $f$. Involving a third arbitrary temperature $T_0$ and heat $Q_0$, the above arguments can be repeated, yielding
\bq
\frac{Q_c}{-Q_0}=f(T_0,T_c), \quad \frac{Q_h}{-Q_0}=f(T_0,T_h)
\eq
Defining $\tau(T):=f(T_0,T)$ this implies
\bq
\frac{Q_h}{-Q_c} = \frac{\tau(T_h)}{\tau(T_c)}
\eq
The function $\tau(T)$ amounts to a {\it re-scaling} of the temperatures (in fact, it corresponds to the absolute thermodynamic scale of temperature \cite{fermi}). For convenience we will again use the same notation $T$ for the rescaled temperature $\tau(T)$. Then \eqref{carnote} for a Carnot cycle can be written as
\bq
\frac{Q_h}{-Q_c} = \frac{T_h}{T_c}
\eq
or equivalently
\bq
\label{carnote2}
\frac{Q_h}{T_h} + \frac{Q_c}{T_c} = 0
\eq
Furthermore, \eqref{carnotin} for an arbitrary (not necessarily Carnot) cycle $(Q_h,T_h; Q_c, T_c)$ amounts to
\bq
\label{carnotin2}
\frac{Q_h}{T_h} + \frac{Q_c}{T_c} \leq 0
\eq
\begin{remark}
\label{universality}
The Second Law is a statement about {\it all} thermodynamic systems; not just a particular system as we are used to in systems and control. This allows us to consider thermodynamic processes consisting of a cycle of one engine and a reversed Carnot cycle of second one. In particular, this {\it universality} implies that the re-scaling of the temperature holds for all thermodynamic systems in the same way, leading to a uniform absolute temperature.
\end{remark}
\subsection{Maximal efficiency of the Carnot cycle}
\label{subsec:efficiency}
The {\it efficiency} of a cycle is defined as the performed mechanical work divided by the absorbed heat at high temperature, given as
\bq
\frac{W}{Q_h}= \frac{Q_h + Q_c}{Q_h} = 1 - \frac{-Q_c}{Q_h}
\eq
Hence it follows from \eqref{carnotin} that Carnot cycles enjoy {\it maximal efficiency} among all cycles, and by \eqref{carnote} this efficiency is independent of which Carnot cycle we take (and by Remark \ref{universality} independent of the system we consider). Furthermore, the efficiency of a Carnot cycle is equivalently given by
\bq
1- \frac{T_c}{T_h}= \frac{T_h-T_c} {T_h}
\eq
In particular, this means that in order to increase efficiency it is most advantageous to lower $T_c$. On the other hand, in many applications (such as the classical the steam engine), $T_c$ is just the temperature of the environment, which can{\it not} be controlled.

\subsection{Clausius' inequality}
\label{subsec:clausius}
The fundamental equality \eqref{carnote2} for a Carnot cycle $(Q_h,T_h; Q_c, T_c)$, and the inequality \eqref{carnotin2} for an arbitrary cycle $(Q_h,T_h; Q_c, T_c)$, can be generalized as follows (see \cite{fermi} for more details). Consider a {\it complex} cycle $(Q_1,T_1; \cdots ;Q_n,T_n)$ consisting of $n$ isothermals at temperatures $T_i$ and absorbed heat quantities $Q_i$, $i=1,2, \cdots,n$, interlaced by $n$ adiabatics. Such a complex cycle $(Q_1,T_1; \cdots ;Q_n,T_n)$ is called {\it heat supply reversible} if also $(-Q_n,T_n; \cdots; -Q_1,T_1)$ is a feasible cycle (returning to the same state). Now let us consider an auxiliary heat source with temperature $T_0$, and $n$ Carnot cycles $(Q_i,T_i; Q_{i,0}, T_0)$,  operating between the temperatures $T_i$ and $T_0$, $i=1,\cdots,n$. According to \eqref{carnote2}
\bq
\label{tempi}
-Q_{i,0}=\frac{T_0}{T_i}Q_i
\eq
Now consider the total transformation consisting of the complex cycle $(Q_1,T_1; \cdots, Q_n,T_n)$ together with the heat supply reversed Carnot cycles $(-Q_{i,0},T_0; -Q_i, T_i), i=1,\cdots,n$. Then the net exchange of heat with each of the sources with temperatures $T_1,\cdots,T_n$ is zero, while the auxiliary source at temperature $T_0$, in view of \eqref{tempi}, absorbs a total heat
\bq
\sum_{i=1}^{n}-Q_{i,0}= T_0 \sum_{i=1}^{n}\frac{Q_i}{T_i}
\eq
But then by Kelvin's formulation of the Second Law this quantity should be less than or equal to zero, or equivalently,
\bq
\label{clausiusin}
\sum_{i=1}^{n}\frac{Q_i}{T_i} \leq 0
\eq
Furthermore, if the complex cycle $(Q_1,T_1; \cdots ; Q_n,T_n)$ is heat supply reversible, we analogously prove the opposite inequality, thus yielding
\bq
\label{clausiuse}
\sum_{i=1}^{n}\frac{Q_i}{T_i} = 0
\eq
A slight extension (approximating continuous heat flow time-functions $q(\cdot)$ by step functions with step values $Q_1,\cdots,Q_n$) then yields the celebrated {\it Clausius inequality}
\bq
\label{Clausiusin}
\oint \frac{q(t)}{T(t)} dt \leq 0
\eq
for all cyclic processes $\big(q(\cdot),T(\cdot)\big)$, with equality 
\bq
\label{Clausiuse}
\oint \frac{q(t)}{T(t)} dt = 0
\eq
holding for heat supply reversible cyclic processes. 

\section{From Clausius' Inequality to Entropy}
\label{sec:entropy}
From the point of view of cyclo-dissipativity theory (see the sidebar "\nameref{sb:cyclo}") the Clausius inequality \eqref{Clausiusin} is exactly the same as {\it cyclo-dissipativity} of any thermodynamic system with respect to the supply rate $-\frac{q}{T}$, where $q$ is the heat flow (heat per second) into the thermodynamic system, and $T$ is the temperature of the system. 
Thus assuming reachability from and controllability from some ground state $x^*$ this means, see Theorem \ref{th:cyclo}, that there exists a {\it storage function} $F$ such that $F(x(t_2)) \leq F(x(t_1)) + \int_{t_1}^{t_2} -\frac{q(t)}{T(t)} dt$.
Hence the function of the state $S:=-F$ satisfies
\bq
\label{Clausiusstate}
S(x(t_2))-S(x(t_1)) \geq \int_{t_1}^{t_2} \frac{q(t)}{T(t)} dt,
\eq
and equivalently (assuming $S$ to be differentiable) its differential version
\bq
\label{Clausiusstate1}
\frac{d}{dt}S \geq \frac{q}{T}
\eq
The function $S$ was called by Clausius '{\it entropy}', from the Greek word $\tau \!\rho o \pi \eta$ for 'transformation'. 
Note that the dissipativity formulation \eqref{Clausiusstate} of the Second Law already appears in \cite{willems1972}, see also \cite{haddad}, although in \cite{willems1972} it is assumed that $F$ is bounded from below, and thus $S$ is bounded from above; corresponding to dissipativity instead of cyclo-dissipativity.

From the point of view of cyclo-dissipativity theory, the storage function $F$ need not be unique. In order to guarantee uniqueness of $F$, and therefore of the entropy $S$ (very desirable from a physics point of view), we may exploit Proposition \ref{prop:revers}. Indeed, once we additionally assume that, given some ground state, for every thermodynamic state there exists a cyclic transformation through this state and the ground state satisfying
\bq
\label{assumption}
\oint \frac{q(t)}{T(t)} =0,
\eq
then by Proposition \ref{prop:revers} the entropy $S$ is indeed unique (up to a constant). The uniqueness of $S$ is, explicitly or implicitly, always assumed in expositions of classical thermodynamics. Cyclo-dissipativity theory thus provides a solid basis for this assumption.

According to \cite{kondepudi} Clausius interpreted the term $\frac{q}{T}$ as the part of the infinitesimal transformation $\frac{d}{dt}S$ that is {\it compensated} by the opposite rate of change $-\frac{q}{T}$ of the {\it entropy of the surrounding}; that is, of the reservoir supplying the heat to the thermodynamic system. 
The remaining part 
\bq
\sigma:= \frac{d}{dt}S - \frac{q}{T} \geq 0
\eq
was called by Clausius the 'uncompensated transformation' ('unkompensierte Verwandlung' in German) \cite{kondepudi}. From this viewpoint the Second Law expresses the fact that the uncompensated transformation is always nonnegative.

If $\sigma$ is non-zero then $\oint \frac{q(t)}{T(t)} < 0$ along any corresponding cyclic path. This implies that at constant temperature the thermodynamic system expels a positive amount of heat to its surrounding; due to an irreversible conversion of, e.g., mechanical, energy into heat. 
The quantity $\sigma \geq 0$ is also called the {\it irreversible entropy production}, and will be the starting point for a broader discussion of irreversible thermodynamics in the section "\nameref{sec:irreversible}". 

\begin{remark}
Interestingly, by a direct application of dissipativity theory, see again the sidebar "\nameref{sb:cyclo}" and in particular Theorem \ref{tm:diss}, it follows that the system is dissipative with respect to the supply rate $- \frac{q}{T}$ if and only if for all states $x$
\bq
\label{dis}
F_a(x)= \sup \int \frac{q(t)}{T(t)}dt < \infty, 
\eq
where the supremum is taken over all heat flow functions $q(\cdot)$ and corresponding temperature profiles $T(\cdot)$ resulting from $x(0)=x$. In fact, if \eqref{dis} holds then $-F_a$ is {\it maximal} among all non-positive functions satisfying \eqref{Clausiusstate}.
\end{remark}

\subsection{Back to the Carnot cycle}
\begin{figure}[h]
\centering
\begin{tikzpicture}

\tikzstyle{bullet} = [circle, fill,minimum size=2pt, inner sep=0pt, outer sep=0pt];

\begin{scope}[shift={(0,0)}]

\draw  (1,2) -- (1,0) -- (4,0);
\node at (4,-0.25) {$S \rightarrow$};
\node at (0.5,1.9) {$T$};
\node at (0.5,2.25) {$\uparrow$};

\node[bullet] (v1) at (1.75,1.75) {};
\node[bullet] (v2) at (1.75,0.75) {};
\node[bullet] (v3) at (3.25,0.75) {};
\node[bullet] (v4) at (3.25,1.75) {};

\begin{scope}[decoration={
    markings,
    mark=at position 0.55 with {\arrow[scale=1.2]{angle 60}}}
    ]
\draw[postaction={decorate}] (v2) -- (v1);
\draw[postaction={decorate}] (v3) -- (v2);
\draw[postaction={decorate}] (v4) -- (v3);
\draw[postaction={decorate}] (v1) -- (v4);
\end{scope}

\node at (2.5,2) {\footnotesize$Q_h$};
\node at (2.5,0.4556) {\footnotesize$Q_c$};

\draw[dashed,gray] (v2) -- (1,0.75);
\draw[dashed,gray] (v1) -- (1,1.75);

\node at (0.825,0.75) {\footnotesize$T_c$};
\node at (0.825,1.75) {\footnotesize$T_h$};

\end{scope}

\begin{scope}[shift={(5,0)}]

\draw  (1,2) -- (1,0) -- (4,0);
\node at (4,-0.25) {$V \rightarrow$};
\node at (0.75,1.9) {$P$};
\node at (0.75,2.25) {$\uparrow$};

\node[bullet] (v4) at (1.5,1.75) {};
\node[bullet] (v3) at (2,0.75) {};
\node[bullet] (v2) at (3.25,0.5) {};
\node[bullet] (v1) at (2.75,1.5) {};

\begin{scope}[decoration={
    markings,
    mark=at position 0.6 with {\arrow[scale=1.2]{angle 60}}}
    ] 
\draw[-,postaction={decorate}] (v1) to [out=-75,in=125] (v2);
\draw[-,postaction={decorate}] (v2) to [out=175,in=-15] (v3);
\draw[-,postaction={decorate}] (v3) to [out=125,in=-75] (v4);
\draw[-,postaction={decorate}] (v4) to [out=-15,in=175] (v1);
\end{scope}

\node at (2.2016,1.9364) {\footnotesize$Q_h, T_h$};
\node at (2.6036,0.3558) {\footnotesize$Q_c, T_c$};

\end{scope}

\draw[thick,-{angle 60}] (4,1.25) .. controls (4.75,1.75) and (5.5,1.25) .. (5.5,1.25);

\end{tikzpicture}

\caption{Carnot cycle in $(V,P)$ and $(S,T)$ diagram}
\label{fig:carnotST}
\end{figure}
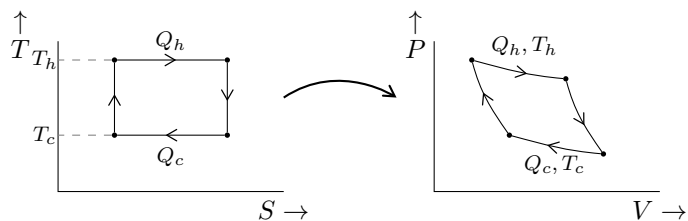

The introduction of the entropy $S$ sheds new light on the Carnot cycle. Indeed, the closed curve in the $(V,P)$ diagram consisting of two isothermals interlaced with two adiabatics, corresponds in the $(S,T)$ diagram to a very simple rectangular curve; see Figure \ref{fig:carnotST}. Furthermore, the resulting map from the $(S,T)$ to the $(V,P)$ diagram is {\it area-preserving}: the area within the rectangular closed curve in the $(S,T)$ diagram (the net amount of absorbed heat $Q_h + Q_c$) is by the First Law equal to the area in the $(V,P)$ diagram circumscribed by the Carnot cycle (the amount of work $W$ done by the system on the surrounding).

\begin{figure}[h]
\centering
\usetikzlibrary{decorations.markings}
\usetikzlibrary{arrows.meta}

\begin{tikzpicture}

\tikzstyle{bullet} = [circle, fill,minimum size=2pt, inner sep=0pt, outer sep=0pt];

\draw  (1,3) -- (1,0) -- (5,0);
\node at (5,-.25) {$V \rightarrow$};
\node at (0.75,2.90) {$P$};
\node at (0.75,3.25) {$\uparrow$};

\node[bullet] (v4) at (1.75,2.25) {};
\node[bullet] (v3) at (2.25,1) {};
\node[bullet] (v2) at (4.2,0.55) {};
\node[bullet] (v1) at (3.7,1.8) {};

\begin{scope}[decoration={
    markings,
    mark=at position 0.55 with {\arrow[scale=1.2]{angle 60}}}
    ] 
\draw[-,postaction={decorate}] (v1) to [out=-80,in=120] (v2);
\draw[-,postaction={decorate}] (v2) to [out=175,in=-22] (v3);
\draw[-,postaction={decorate}] (v3) to [out=120,in=-80] (v4);
\end{scope}

\begin{scope}[decoration={
    markings,
    mark=at position 0.4 with {\arrow[very thick,scale=1]{Bar}}}
    ] 
\draw[-,postaction={decorate}] (v4) to [out=-22,in=175] (v1);
\end{scope}

\node at (3.2655,2.1466) {\footnotesize$-Q_c, T_h$};
\node at (2.1599,2.477) {\footnotesize$Q_h + Q_c$};
\node at (3.1982,0.4487) {\footnotesize$Q_c, T_c$};

\end{tikzpicture}
\caption{Carnot cycle re-interpreted}
\label{fig:carnotre}
\end{figure}
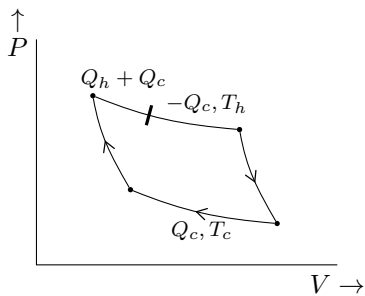

Finally, let us close this section with a simple, but intriguing, re-interpretation of the fundamental equality \eqref{carnote2} for a Carnot cycle $(Q_h,T_h; Q_c, T_c)$; see also \cite{chinees} for similar reasoning in the original work of Clausius. Rewrite \eqref{carnote2} as
\bq
0 = \frac{Q_h}{T_h} + \frac{Q_c}{T_c} = \frac{Q_h +Q_c}{T_h} + Q_c \big(\frac{1}{T_c} - \frac{1}{T_h}\big),
\eq
or equivalently
\bq
\label{alternativecarnot}
\frac{Q_h+Q_c}{T_h} = - Q_c \big(\frac{1}{T_c} - \frac{1}{T_h}\big)
\eq
Note that $-Q_c$ is the net amount of heat flowing from the hot heat source to the cold one, while $\frac{1}{T_c} - \frac{1}{T_h}$ is known as the {\it thermodynamic force} between the two heat sources; see the section \nameref{sec:irreversible}. Furthermore, 
$Q_h+Q_c$ is the total heat supplied by the two sources to the thermodynamic system (which by the First Law equals the work $W$ performed by the system on its surrounding), and corresponds to the left part of the upper isothermal in the Carnot cycle in Figure \ref{fig:carnotre}. Hence \eqref{alternativecarnot} means that the entropy increase $\frac{Q_h+Q_c}{T_h}$ due to the supplied heat during the left part of the upper isothermal is equal to the change in entropy due to a direct heat flow $-Q_c$ from hot to cold (the falling caloric flow in the original interpretation of Carnot; see the sidebar "\nameref{sb:caloric}").

Apart from the Carnot cycle there are other interesting cycles as well. An extensive discussion of them is given in  \cite{zemansky}. Like the Carnot cycle originates from the steam engine, most of these alternative cycles derive from the operation of other types of engines. For example, the Otto cycle is described in the $(V,P)$ diagram by two adiabatics and two isochores (constant volume), the Diesel cycle by two isobars (constant pressure) and two adiabatics, and the Stirling cycle by two isothermals and two isochores.

\subsection{Kelvin and the limits to energy conversion}
\label{subsec:kelvin}

We already noted in the beginning of this section that Kelvin's formulation of the Second Law implies that for each constant temperature $T$ the thermodynamic system is {\it cyclo-passive} with respect to the supply rate $-Pu_V$ corresponding to mechanical work. 
By the main theorem of cyclo-dissipativity, Theorem \ref{th:cyclo}, this means that for each constant temperature $T$ there exists a function $F_T$ of the state $x$ of the thermodynamic system satisfying
\bq
\label{helmholtz}
F_T(x(t_2)) - F_T(x(t_1))  \leq  \int_{t_1}^{t_2}-P(t)u_V(t) dt
\eq
This expresses the property that no thermodynamic system can convert thermal energy from a heat source of constant temperature into mechanical work while returning to its original state. 
From a general system-theoretic point of view this leads to the following question. Consider a system 
with two 'ports' $(u_1,y_1)$ and $(u_2,y_2)$, which is cyclo-passive with respect to the supply rate $s(u_1,u_2,y_1,y_2)=y_1^\top u_1 + y_2^\top u_2$, i.e., there exists a state function $F$ such that
\bq
\frac{d}{dt}F \leq y_1^\top u_1 + y_2^\top u_2
\eq
Under which conditions is it {\it not} possible to transform energy at port $1$ into energy at port $2$ while keeping $y_1$ constant? Or said differently, what is so special about thermodynamic systems, and are there any systems other than thermodynamic systems that can{\it not} transform energy from one port into the other while keeping the output at the first port constant? 
The answer to this question is: Yes, there are quite a few other systems sharing this property. 
In fact, sufficient (and often necessary) conditions on the system for this to happen are discussed in \cite{jeltsema-limits}, together with a number of examples from different areas (synchronous machine, DC-machine, capacitor microphone, $\cdots$). 

Returning to the original dissipation inequality \eqref{helmholtz}, we finally note that even if the storage function $F_T$ for each $T$ is unique, it is so only {\it up to a constant}, and this constant may arbitrarily depend on $T$. This is intimately related to the fact that Kelvin's formulation of the Second Law is {\it stronger} than just cyclo-passivity for every $T$. In fact, the {\it canonical} choice of $F_T$ for every $T$ will be provided by the {\it Helmholtz free energy}; one of the thermodynamic potentials to be discussed in section "\nameref{sec:gibbs}".

\section{Irreversible Thermodynamics}
\label{sec:irreversible}

Recall the dissipation inequality $\frac{d}{dt}S \geq \frac{q}{T}$, equivalently re-written as
\bq
\frac{d}{dt}S = \frac{q}{T} + \sigma; \quad \sigma \geq 0,
\eq
where $\sigma := \frac{d}{dt}S - \frac{q}{T}$ is the {\it irreversible entropy production} \cite{kondepudi} ('uncompensated transformation' in the terminology of Clausius). 
Irreversible thermodynamics is concerned with the dynamics of thermodynamic systems in case $\sigma$ is different from zero, implying an autonomous (independent from external heat flow) increase of the entropy $S$. Sometimes it is also referred to as {\it non-equilibrium thermodynamics}, because of resulting from (internal) non-equilibrium conditions.

The perhaps simplest example of irreversible dynamics and irreversible entropy production is offered by the {\it heat exchanger}. Consider two heat compartments, having temperatures $T_h$ and $T_c$ ('hot' and 'cold'), connected by a  heat-conducting wall. In the absence of the conducting wall, these are two separate systems with entropies $S_h$ and $S_c$, each satisfying
\bq
\frac{d}{dt}S_h= \frac{q_h}{T_h}, \quad \frac{d}{dt}S_c= \frac{q_c}{T_c}
\eq
Due to the conducting wall there is a heat flow $q$ from the hot to the cold compartment, while in view of the First Law $q=-q_h=q_c$. Hence the total entropy $S:=S_h + S_c$ satisfies
\bq
\frac{d}{dt}S = -\frac{q}{T_h} + \frac{q}{T_c}= \big(\frac{1}{T_c} - \frac{1}{T_h} \big)q ,
\eq
where, using Fourier's law for heat conduction, $q= \lambda (T_h - T_c)$ for some positive constant $\lambda$. This yields the following expression for the irreversible entropy production $\sigma = \frac{d}{dt}S$ (note that there is no external heat flow)
\bq
\sigma = \big(\frac{1}{T_c} - \frac{1}{T_h} \big) \lambda (T_h - T_c)= \lambda \frac{\big(T_h - T_c \big)^2}{T_hT_c} \geq 0
\eq
A second simple example refers to the adiabatics in the classical Carnot cycle of a gas in a volume controlled by a piston. If the pressure $P$ of the gas is not considered to be the direct control variable, but instead one distinguishes between a pressure $P_{\mathrm{piston}}$ on the piston, and a pressure $P_{\mathrm{gas}}$ of the gas, then during the expansion phase $P_{\mathrm{gas}} \geq P_{\mathrm{piston}}$, leading to the irreversible entropy production $\sigma$ given by
\bq
\sigma= \frac{P_{\mathrm{gas}} - P_{\mathrm{piston}}}{T} J ,
\eq
where $J$ is the rate of volume change. Typically, $J$ will be positively proportional to the difference $P_{\mathrm{gas}} - P_{\mathrm{piston}}$, implying that, indeed, $\sigma\geq 0$.

\subsection{Irreversibility of chemical reaction networks}
\label{subsec:chemical}
A much more involved example of irreversible thermodynamics is the dynamics of non-isothermal {\it chemical reaction networks}. 
Consider for simplicity an {\it isolated} chemical reaction network (no incoming or outgoing mass flow, and no external heat flow), with $m$ chemical species and $r$ reactions. Also, let us disregard volume and pressure, and model the state of the chemical reaction network by the vector $x \in \mR^m$ of concentrations of the chemical species. Then the dynamics of the concentrations takes the form
\bq
\dot{x}=Nv(x),
\eq
where $N$ is an $m \times r$ matrix, called the {\it stoichiometric} matrix, and $v \in \mathbb{R}^r$ is the vector of reaction fluxes. The stoichiometric matrix $N$, which consists of (positive and negative) integer elements, captures the basic conservation laws of the reactions. Chemical reaction network theory, originating from \cite{HornJackson, Horn} and \cite{Feinberg1}, identifies the edges of the underlying directed graph with the reactions, and the nodes with the $c$ {\it complexes} of the chemical reactions, i.e., all the different left- and right-hand sides of the reactions in the network. This means the stoichiometric matrix $N$ is factorized as $N=ZB$, with $B$ denoting the $c \times r$ {\it incidence} matrix of the graph of complexes, and $Z$ the $m \times c$ {\it complex composition} matrix (a matrix of nonnegative integers), whose $\rho$-th column captures the expression of the $\rho$-th complex in the $m$ chemical species. It is shown in \cite{SIAP} that the dynamics $\dot{x}=Nv(x)$ of a large class of chemical reaction networks (including detailed-balanced mass action kinetics networks) can be written into the compact form
\bq
\label{masterequation1}
\dot{x} = - Z \ML \Exp \left(\frac{1}{RT}Z^\top \mu(x)\right),
\eq
where $\Exp$ is the vector exponential mapping $\Exp (z)=(\exp z_1,\cdots,\exp z_c)^\top$, $R$ is the gas constant, $T$ is the temperature, and $\mu$ is the $m$-dimensional vector of {\it chemical potentials} of the chemical species (for which e.g. in the case of detailed-balanced mass action kinetics an explicit expression is available). 
Furthermore, the matrix $\ML:= B \mathcal{K}B^\top$ in \eqref{masterequation1} defines a {\it weighted Laplacian} matrix for the graph of complexes, with the diagonal elements of the diagonal matrix $\mathcal{K}$ given by the so-called conductances $\kappa_1, \cdots, \kappa_r$ (which are depending on the temperature  $T$ and the reference state).  
We have the following fundamental property \cite{SIAP}
\bq
\label{fund}
\gamma^\top \ML \Exp \, \gamma \geq 0 \mbox{ for all }\gamma \in \mR^c, \quad \gamma^\top \ML \Exp \, \gamma = 0 \mbox{ iff } B^\top \gamma=0
\eq
The entropy $S$ expressed as a function of $x$ and the total energy $E$ satisfies by Gibbs' fundamental relation
\bq
\frac{\partial S}{\partial x}(x,E)= - \frac{\mu}{T}, \quad \frac{\partial S}{\partial E}(x,E)=\frac{1}{T}
\eq
This implies that for an isolated chemical reaction network (no external heat or matter flow; hence constant energy $E$)
\bq
\frac{d}{dt}S = \frac{1}{T} \mu^\top Z \ML \Exp (\frac{Z^\top \mu}{RT}) =: \sigma \geq 0,
\eq
with equality if and only if $B^\top Z^\top \mu=N^\top \mu=0$, i.e., if and only if the {\it chemical affinities} $N^\top \mu$ of the reactions are all zero. Hence the equilibria of the system correspond to states of minimal (i.e., zero) entropy production $\sigma$, in accordance with the theory of irreversible thermodynamics \cite{kondepudi}. By using $-S$ as Lyapunov function, it follows, under the standard assumption that trajectories will not converge to the boundary of the positive orthant $\mathbb{R}^m_+$, that any initial vector of concentrations in the positive orthant will converge to one of these equilibria; see e.g. the exposition in \cite{SIAP, wang, louvain}.

\subsection{Thermodynamic forces and flows}
\label{subsec:forces/flows}
In irreversible thermodynamics it is normally postulated \cite{kondepudi} that the irreversible entropy production can be represented as
\bq
\label{factorization}
\sigma = \sum_{k=1}^s F_kJ_k \geq 0,
\eq
where $F_k$ are the {\it thermodynamic forces} and $J_k$ are the {\it thermodynamic flows} (or fluxes), in such a way that
\bq
\sigma =0 \Leftrightarrow F_k=0, k=1, \cdots, s
\eq
The examples given above do illustrate this postulate. In the heat exchanger example the thermodynamic force is $F=\frac{1}{T_c} - \frac{1}{T_h}$, while the thermodynamic flow is $q= \lambda (T_h - T_c)$. Note that indeed $\sigma=0$ if and only if $F=$. In the piston example, the thermodynamic force is $F= \big(P_{\mathrm{gas}} - P_{\mathrm{piston}}\big)\cdot \frac{1}{T}$ and the flow is $u_V = \mu F$ for some $\mu >0$.

In the case of isolated chemical reaction networks, the vector of thermodynamic forces $F$ is given as $F= \frac{1}{T}N^\top \mu$, the vector of chemical affinities divided by temperature $T$. Furthermore, the vector of thermodynamic flows $J$ is given as 
\bq
\label{eq:flow}
J= \mathcal{K}B^\top \Exp \frac{Z^\top \mu}{RT}
\eq
(which is also equal to the vector of {\it rates of extent} of each reaction). Furthermore, it follows from \eqref{fund} that $\sigma=0$ if and only if $F=0$. 
From a systems and control point of view, the factorization \eqref{factorization} of $\sigma = \frac{d}{dt}S - \frac{q}{T}$ is closely related to the factorization of the differential dissipation inequality \eqref{diss-eq} in (cyclo-)dissipativity theory \cite{willems1972,hillmoylan1980}. For example, in the case of linear systems with quadratic supply rate, the storage function is also quadratic, and the differential dissipation inequality amounts to a Linear Matrix Inequality (LMI), which can be factorized as in \eqref{factorization}. Similarly, in {\it linear} irreversible thermodynamics \cite{kondepudi} it is assumed that the vector $F$ of thermodynamic forces and the vector $J$ of thermodynamic flows are linearly related as
\bq
\label{Onsager}
J=LF, \quad L=L^\top
\eq
These relations, and especially the symmetry of the matrix $L$, are the celebrated {\it Onsager reciprocity relations} \cite{kondepudi}. They lead to the symmetric factorization $\sigma= F^\top J= F^\top LF$. 

Note that in the piston example the thermodynamic flow $J$ {\it is} indeed expressed as $J = \mu F$, and is (trivially) satisfying the Onsager relations. However the heat exchanger is {\it not} of this form. In fact, the thermodynamic flow $q= \lambda (T_h - T_c)$ can{\it not} be expressed as a function of the thermodynamic force $F=\frac{1}{T_c} - \frac{1}{T_h}$ (although $q=0$ if and only if $F=0$). Similarly, chemical reaction networks are {\it not} of this form. The relation between $F$ and $J$ in this case is {\it not} linear, and in most cases $J$ as in \eqref{eq:flow} can{\it not} be expressed as a function of $F$.

\subsection{Cyclo-passive and port-Hamiltonian systems as irreversible thermodynamic systems}
\label{subsec:cyclo-passive}
A large class of systems which, somewhat artificially, can be formulated within an irreversible thermodynamics context are standard {\it input-state-output systems}
\bq
\label{iso}
\begin{array}{rcl}
\dot{x} & = & f(x) + g(x)u, \quad u \in \mR^m, \\[2mm]
y & = & h(x) , \quad y \in \mR^m,
\end{array}
\eq
where $x$ is in some $n$-dimensional state space manifold $\X$, which are assumed to be {\it cyclo-passive}; see the sidebar "\nameref{sb:cyclo}". This means that there exists a state function $H$ such that $\frac{d}{dt}H \leq y^\top u$, i.e.,
\bq
\label{passivity}
\frac{\partial H}{\partial x^\top}(x)f(x)=: -\rho(x) \leq 0, \quad \frac{\partial H}{\partial x^\top}(x)g(x) = h^\top (x)
\eq
Now define an {\it additional state variable} $S$ (interpreted as the {\it entropy} of the system), together with an 'internal energy' $U(S)$. Then consider the {\it total} energy $E(x,S):=H(x) + U(S)$, and extend the system \eqref{iso} to
\bq
\label{isoS}
\begin{array}{rcl}
\dot{x} & = & f(x) + g(x)u \\[2mm]
\dot{S} & = & \frac{\rho(x)}{U'(S)} \\[2mm]
y & = & h(x)
\end{array}
\eq
This extended system satisfies 
\bq
\frac{d}{dt}E= \frac{d}{dt}H + \frac{d}{dt}U= -\rho(x) +y^\top u + U'(S)\frac{\rho(x)}{U'(S)}=y^\top u,
\eq
and thus is cyclo-lossless; satisfying the First Law. Furthermore
\bq
\frac{d}{dt}S  = \frac{\rho(x)}{U'(S)}=: \sigma \geq 0 ,
\eq
corresponding to the Second Law. Note that the choice of the internal energy function $U(S)$ is rather arbitrary, with $U'(S)$ defining an (artificial) {\it temperature}. One possible choice is $U(S)=T_0S$ corresponding to an infinite heat reservoir at constant temperature $T_0$. 

The factorization \eqref{factorization} for cyclo-passive systems becomes most clear if the cyclo-passive system \eqref{iso} can be represented into {\it port-Hamiltonian} form \cite{maschkebordeaux,passivitybook,schjeltsema}
\bq
\label{pH}
\begin{array}{rcl}
\dot{x} & = & \J(x)\frac{\partial H}{\partial x}(x) - g_R(x)R\big( g_R^\top (x)\frac{\partial H}{\partial x}(x)\big)+ g(x)u \\[2mm]
y & = & g^\top(x) \frac{\partial H}{\partial x}(x)
\end{array}
\eq
for some mapping $R$ satisfying $z^\top R(z)\geq 0$ for all vectors $z=g_R^\top (x)\frac{\partial H}{\partial x}(x)$, a skew-symmetric matrix $\J(x)$, and matrices $g(x),g_R(x)$. Here the term $- g_R(x)R\big( g_R^\top (x)\frac{\partial H}{\partial x}(x)\big)$ models {\it energy dissipation} (without taking into account the produced heat). Note that any port-Hamiltonian system is cyclo-passive with storage function $H$, since the requirement $z^\top R(z)\geq 0$ implies $\frac{d}{dt}H \leq y^\top u$. Considering as before an additional entropy variable $S$ and internal energy $U(S)$, the irreversible entropy production $\sigma$ takes the form
\bq
\sigma = \frac{\partial H}{\partial x^\top}(x)g_R(x)R\big(g_R^\top (x)\frac{\partial H}{\partial x}(x)\big)\cdot\frac{1}{U'(S)} \geq 0,
\eq
which is already in factorized form $F^\top J$ with
\bq
F=g_R^\top (x) \frac{\partial H}{\partial x }(x)\cdot\frac{1}{U'(S)}, \quad J= R\big(g_R^\top (x)\frac{\partial H}{\partial x}(x)\big)
\eq
In fact, if $R$ is a symmetric linear mapping then the Onsager reciprocity relations \eqref{Onsager} are satisfied.

\section{Gibbs and the Thermodynamic Phase Space}
\label{sec:gibbs}
This section marks the transition to the {\it geometrization} of classical thermodynamics. In particular, Gibbs' fundamental relation (between the extensive and intensive variables) defines the constitutive relations (state properties) of thermodynamic systems. This naturally leads to a contact-geometric formulation, initiated by Gibbs and explicitly stated by Hermann \cite{hermann}, and entails a 'paradigm shift' towards non-minimal systems modeling by the introduction of the 'thermodynamic phase space'.

\subsection{Gibbs' fundamental thermodynamic relation and thermodynamic potentials}
\label{subsec:gibbs}
Consider a simple thermodynamic system, with variables $V,P,T$. The {\it equation of state} is an equation  $f(V,P,T)=0$ for some scalar function $f$; see subsection \nameref{subsec:basic}. Throughout we assume that the thus defined set of states of the thermodynamic system is a $2$-dimensional {\it submanifold} $M$ of $\mR^3$ (strictly speaking of $\mR^2 \times \mR^+$, since $T\geq 0$). 

Using the First and Second Law we defined functions $E: M \to \mR$ ({\it energy}) and $S: M \to \mR$ ({\it entropy}). Here $E$ is unique up to a constant, while also $S$ is unique up to a constant, under the additional assumption, cf. Proposition \ref{prop:revers}, that given some ground state there exists for any state a cyclic path through this state and the ground state satisfying, cf. \eqref{assumption},
\bq
\label{assumption1}
\oint \frac{q(t)}{T(t)} =0
\eq
This will be a standing assumption throughout.
Then we may equally well represent the set of states $M \subset \mR^3$ by the $2$-dimensional submanifold $L \subset \mR^5$ given as
\bq
L:= \{(E,S,V,T,P) \mid f(V,P,T)=0, E=E(V,P,T), S= E(V,P,T) \}
\eq
Note that, with some abuse of notation, we have introduced here the extra {\it variables} $E,S$, denoted by the same letters as used for the {\it functions} defined before. Then, under reasonable assumptions, we can parametrize $L$ by the extensive variables $S$ and $V$, and consider the so-called {\it energy representation} of the submanifold $L \subset \mR^5$ given as
\bq
L:= \{(E,S,V,T,P) \mid E=E(S,V), T= T(S,V), P= P(S,V) \}
\eq
for some functions $E(S,V), T(S,V), P(S,V)$. Thus, the extensive variable $E$, as well as the two intensive variables $T$ and $P$ are expressed as functions of the remaining extensive variables $S,V$, which are serving as coordinates for $L$. The space $\mR^5$ consisting of all the variables $E,S,V,T,P$ is called the {\it thermodynamic phase space}.

Now let us exploit once more the First and Second Law. By the First Law $\frac{d}{dt}E=-P\frac{d}{dt}V +q$. Furthermore, by the assumption \eqref{assumption1} there exists for any state a path through this state and the ground state such that
\bq
 \frac{d}{dt}E=-P\frac{d}{dt}V + T\frac{d}{dt}S
\eq
This implies that the {\it Gibbs' one-form} on $\mR^5$ defined as
\bq
\label{Gibbs}
dE - TdS +PdV, \quad \mbox{ Gibbs' one-form }
\eq
is {\it zero} restricted to $L$. This is called {\it Gibbs' fundamental thermodynamic relation}. It implies that the submanifold $L$ is actually given as
\bq
L:= \{(E,S,V,T,P) \mid E=E(S,V), T= \frac{\partial E}{\partial S}(S,V), -P= \frac{\partial E}{\partial V}(S,V) \}
\eq
Thus $L$ is completely described by the energy function $E(S,V)$, whence the name {\it energy representation}.

On the other hand, the submanifold $L$ may equally well be parametrizable by, e.g., the variables $T,V$. Define the {\it partial Legendre transform} of $E(S,V)$ with respect to $S$ as
\bq
A(T,V) :=  E(V,S) -TS, \quad T= \frac{\partial E}{\partial S}(S,V),
\eq
where $S$ is solved from $T= \frac{\partial E}{\partial S}(S,V)$ as a function of $(T,V)$. This means that $L$ is also described as
\bq
L:= \{(E,S,V,T,P) \mid E=A(T,V) - T \frac{\partial A}{\partial T}(T,V), S= - \frac{\partial A}{\partial T}(T,V), -P= \frac{\partial A}{\partial V}(T,V) \}
\eq
The function $A(V,T)$ is called the {\it Helmholtz free energy}, and is one of the thermodynamic potentials, derived from $E(S,V)$, to describe the submanifold $L$.
For example, in case of an ideal gas \cite{fermi, kondepudi}
\bq
E(S, V)= \frac{C_Ve^{\frac{S}{C_V}}}{Ve^{\frac{R}{C_V}}},
\eq
where $C_V$ denotes the heat capacity (at constant volume), and $R$ is the universal gas constant. Partial Legendre transform $E(S, V)$ with respect to $S$ yields the Helmholtz free energy $A(T,V)$ given as \cite{fermi}
\bq
A(T,V) = C_VT + W - T\big(C_V \ln T + R \ln V +a\big),
\eq
for constants $a$ (the entropy constant of the gas) and $W$ (an integration constant). 

Apart from the energy $E(S,V)$ and the Helmholtz free energy $A(T,V)$ there are two more thermodynamic potentials that can be obtained from $E(S,V)$ by partial Legendre transform:
\bq
\begin{array}{rclll}
H(S,P) & = & E(S,V)+PV, \quad & \mbox{enthalpy},& \mbox{coordinates } P,S \\[4mm]
G(T,P) &= & H(S,P) -TS, \quad & \mbox{Gibbs' free energy}, & \mbox{coordinates } P,T
\end{array}
\eq
Despite all these different ways to parametrize $L$ by two coordinates, corresponding to different thermodynamic potentials as described above, the situation is {\it very simple} from a {\it geometric} point of view: there is just {\it one} $2$-dimensional submanifold $L$ describing the set of states of the thermodynamic system, which is such that the Gibbs' form $dE - TdS +PdV$ is zero restricted to it. The appropriate geometric setting for all this is {\it contact geometry}, as already alluded to in the quote by Arnold in the Introduction, and first explored within thermodynamics in e.g. \cite{hermann,mrugala1,mrugala2}; see \cite{bravetti} for a survey on recent developments. The use of contact geometry for thermodynamic systems from a control point of view was initiated in \cite{eberard}. The basic notions of contact geometry are discussed in the sidebar "\nameref{sb:contact}". 

Apart from the above options to parametrize the submanifold $L$ by different sets of coordinates, corresponding to different thermodynamic potentials derivable from the energy $E(S,V)$, there is an {\it alternative}, but similar, way of describing the set $L$ of thermodynamic states. This is to start with the expression of the {\it entropy} as a function $S(E,V)$ of the volume and the energy, and leads to the {\it entropy representation}. This alternative option can be also motivated from a modeling point of view. Namely, in many situations (e.g., within chemical engineering) thermodynamic systems are formulated by first listing the {\it balance laws} for the extensive variable $V$, and the mole numbers $N_k$ of the chemical species, as well as the energy $E$. Then the entropy is sought to be expressed as a function of these extensive variables. For a simple thermodynamic system this leads to the representation of the submanifold $L \subset \mR^5$ given as
\bq
L:= \{(E,S,V,T,P) \mid S=S(E,V), \frac{1}{T}= \frac{\partial S}{\partial E}(E,V), \frac{P}{T}= \frac{\partial S}{\partial V}(E,V) \}
\eq
Starting from this entropy representation, one may then define, as in the energy representation, other thermodynamic potentials obtained by partial Legendre transform of $S(E,V)$.
Also for the entropy representation the geometric point of view is the simplest. It corresponds to the one-form
\bq
dS -\frac{1}{T}dE - \frac{P}{T}dV
\eq
being zero on $L$. This one-form is obtained from the Gibbs one-form $dE - TdS +PdV$ by dividing by $-T$, and leads to a different, although very similar, contact-geometric description.

Gibbs' fundamental thermodynamical relation is immediately extended to more general situations. For instance, see \cite{openchemical} for a careful derivation, in the case of multiple chemical species with mole numbers $N_1, \cdots, N_m$ and chemical potentials $\mu_1, \cdots, \mu_m$, it amounts to the extended Gibbs one-form
\bq
dE - TdS +PdV - \sum_{k=1}^m  \mu_kdN_k
\eq
to be zero on a submanifold $L$ of the form (in energy representation)
\bq
\begin{array}{rl}
L & := \{(E,S,V,N_1, \cdots, N_k,T,P,\mu_1,\cdots, \mu_k) \mid E=E(S,V,N_1, \cdots,N_k), \\[2mm]
& \; T = T(S,V,N_1, \cdots,N_k), -P = P(S,V,N_1, \cdots,N_k), \mu_i=\mu_i(S,V,N_1, \cdots,N_k) \},
\end{array}
\eq
which implies that
\bq
\begin{array}{rl}
L  & :=  \{(E,S,V,N_1, \cdots, N_k,T,P,\mu_1,\cdots, \mu_k) \mid E=E(S,V,N_1, \cdots,N_k),\\[2mm]
& \qquad T =  \frac{\partial E}{\partial S}(S,V,N_1, \cdots,N_k),-P= \frac{\partial E}{\partial V}(S,V,N_1, \cdots,N_k),\\[2mm]
& \qquad \mu_i= \frac{\partial E}{\partial N_i}(S,V,N_1, \cdots,N_k), \; i=1, \cdots, n \}
\end{array}
\eq
with $E(S,V,N_1, \cdots,N_k)$ the energy function. By partial Legendre transform of $E(S,V,N_1, \cdots,N_k)$ one obtains thermodynamic potentials  corresponding to other parametrizations of $L$. 
Similarly, by expressing the entropy as $S=S(E,V,N_1, \cdots,N_k)$ one obtains the {\it entropy representation} of $L$.

\subsection{A paradigm shift in systems modeling}
\label{subsec:paradigm}

Let us reflect on what we have seen so far in this section. The state space of a simple thermodynamic system is described by a $2$-dimensional submanifold of the $3$-dimensional space of macroscopic quantities $V,P,T$; one extensive, and two intensive. Then, based on the First and Second Law of thermodynamics, two extra extensive variables $E,S$ are introduced, and the state space is equivalently described as a $2$-dimensional submanifold $L$ of $\mR^5$; the space of the three extensive variables $E,S,V$, and the two intensive variables $T,P$. The submanifold $L$ defines the {\it constitutive relations} of the thermodynamic system, that is, the {\it state properties} of the system. Characterizing property of $L$ is that it is a maximal submanifold restricted to which the Gibbs' one-form is zero. Such manifolds are called {\it Legendre submanifolds}; see the sidebar "\nameref{sb:contact}". Furthermore, any such $L$ defines possible constitutive relations. For example, the Legendre submanifold $L$ corresponding to an {\it ideal gas} is different from the Legendre submanifold $L$ corresponding to a {\it Van der Waals gas} \cite{fermi}.
Hence, in general Gibbs' fundamental thermodynamical relation corresponds to the {\it constitutive relations} of the thermodynamic system. The thermodynamic phase space is the total space of all involved variables (extensive and intensive; one more extensive variable than intensive), and care should be taken to regard this as the state space. Instead, the {\it minimal} state space of the thermodynamic system is the Legendre submanifold $L$ of the thermodynamic phase space.

Furthermore, note that no {\it dynamics} is yet defined. The First and Second Law impose {\it constraints} on any possible dynamics. Furthermore, they lead to the definition of the extensive variables energy and entropy, and their combination implies Gibbs' fundamental relation, characterizing all the minimal state spaces (i.e., the Legendre submanifolds $L$).

A more or less appropriate analogy outside the usual thermodynamic realm is the following. Consider a capacitor, with charge $Q$, and $V$ the voltage across the capacitor. Then the constitutive relations of the capacitor are specified by an energy function $E(Q)$, yielding the $1$-dimensional submanifold
\bq
\widetilde{N}= \{ (Q,V) \mid V= \frac{d E}{d Q}(Q) \} \subset \mR^2 
\eq
$Q$ can be considered to be an extensive variable, and $V$ an intensive variable. By introducing the energy $E$ as an {\it extra} extensive variable, this leads to the equivalent description of the capacitor by the $1$-dimensional submanifold
\bq
N= \{ (E,Q,V) \mid E=E(Q),V= \frac{d E}{d Q}(Q) \} \subset \mR^3,
\eq
in the {\it extended} space of two variables $E,Q$, and intensive variable $V$. The submanifold $N$ defines a maximal submanifold of $\mR^3$ restricted to which the one-form $dE-VdQ$ is zero (i.e., a Legendre submanifold); analogously to the Legendre submanifold $L$ of a simple thermodynamic system. Thus for a capacitor the 'thermodynamic phase space' is $\{ (E,Q,V) \in \mR^3 \}$, while the constitutive relations are defined by $N$, or equivalently by the energy function $E(Q)$.

\section{From Thermodynamic Phase Space to Hamiltonian Dynamics}
\label{sec:hamiltonian}
Thermodynamics as discussed so far is basically thermo{\it statics}; dealing with the constitutive relations of the system, ultimately characterized by the Legendre submanifold $L$ of the thermodynamic phase space of all extensive variables (including energy and entropy) and intensive variables (Gibbs' fundamental thermodynamic relation). This implies that any dynamics should be such that the constitutive relations are respected; i.e., any dynamics defined on the thermodynamic phase space should leave the Legendre submanifold characterizing the constitutive properties {\it invariant}. This can be formulated using contact geometry and the notion of a {\it contact vector field} on a contact manifold; see the sidebar "\nameref{sb:contact}". However, in the next subsection we will immediately take one more abstraction step, which at the same time will {\it resolve} some problems in the contact-geometric formulation of thermodynamic systems, and also {\it simplify} the representations and computations. Eventually, it will be also crucial in the definition of {\it thermodynamic ports} in the section "\nameref{sec:port-thermodynamic}". This abstraction step is the step from {\it contact geometry} to {\it homogeneous symplectic geometry}.

\subsection{From contact to homogeneous symplectic geometry}
The contact-geometric view on thermodynamics, despite being directly motivated by Gibbs' fundamental relation, has three shortcomings:\\
(1) Switching from the energy representation $E=E(S,V)$ to the entropy representation $S=S(E,V)$ corresponds to dividing the Gibbs form $dE - TdS +PdV$ by $-T$, leading to the {\it new} contact form
\bq
dS - \frac{1}{T}dE - \frac{P}{T}dV
\eq
with {\it new} intensive variables $\frac{1}{T}, \frac{P}{T}$.
Obviously $L$ is a Legendre submanifold for this new contact form as well, but strictly speaking it leads to a {\it different} contact-geometric description.\\
(2) In general, the contact-geometric approach does not make a clear distinction between extensive and intensive variables; given a contact form $\theta$ there are many Darboux coordinates as in \eqref{darboux} in the sidebar "\nameref{sb:contact}".\\
(3) Computations in contact geometry tend to be involved; especially when it comes to {\it dynamics}, as will be the topic of the next section. 

The way to solve these problems is to {\it extend} contact manifolds by one extra dimension to {\it symplectic} manifolds, in fact cotangent bundles, with an added {\it homogeneity} structure. For a simple thermodynamic system with extensive variables $E,S,V$ and intensive variables $T,-P$, this amounts to replacing the intensive variables $T,-P$ (in the energy representation) by their {\it homogeneous coordinates} $p_E,p_S,p_V$ with $p_E \neq 0$, i.e., 
\bq
T= \frac{p_S}{-p_E}, \; -P= \frac{p_V}{-p_E},
\eq
and thus to express the intensive variables $\frac{1}{T}, \frac{P}{T}$ in the entropy representation as
\bq
\frac{1}{T} = \frac{p_E}{-p_S}, \; \frac{P}{T}= \frac{p_V}{-p_S}
\eq
In this way, the {\it two} Gibbs' one-forms $dE - TdS +PdV$ and $dS -\frac{1}{T}dE - \frac{P}{T}dV$ are replaced by a {\it single} symmetric expression, namely by the Liouville one-form
\bq
p_EdE + p_SdS + p_VdV,
\eq
being the {\it canonical one-form} on the cotangent bundle $T^*\mR^3$, with $\mR^3$ the space of extensive variables $E,S,V$. By the definition of homogeneous coordinates the vector $(p_E,p_S,p_V)$ is different from the $0$-vector. Hence the space $\{(E,S,V,p_E,p_S,p_V) \}$ is the cotangent bundle $T^*\mR^3$ {\it minus} its zero section.
Using homogeneous coordinates the $2$-dimensional Legendre submanifold $L$ is now replaced by the $3$-dimensional submanifold $\cL \subset T^*\mR^3$, given as
\bq
\cL=\{(E,S,V,p_E,p_S,p_V) \mid (E,S,V, \frac{p_S}{-p_E}, \frac{p_V}{-p_E}) \in L \; (p_E,p_S,p_V) \neq 0  \}
\eq
It turns out that $\cL$ is a {\it Lagrangian submanifold}, which is moreover {\it homogeneous}, in the sense that whenever 
$(E,S,V,p_E,p_S,p_V) \in \cL$ then also $(E,S,V,\lambda p_E,\lambda p_S, \lambda p_V) \in \cL$, for any non-zero $\lambda \in \mR$. Such homogeneous Lagrangian submanifolds are fully characterized as maximal manifolds restricted to which the canonical one-form $p_EdE + p_SdS + p_VdV$ is zero; see the sidebar "\nameref{sb:homogeneous}".

In this way, the first two disadvantages of the contact geometry formulation (difference between energy and entropy representation, and the lack of clear distinction between extensive and intensive variables) are resolved. As explained in the sidebar "\nameref{sb:homogeneous}" this symplectization of contact manifolds (by adding one extra dimension to the space of intensive variables) has clear computational advantages as well. In fact, all computations become standard operations in Hamiltonian dynamics. In the words of Arnold \cite{arnold-contact}: one is advised to calculate symplectically (but to think rather in terms of contact geometry). The symplectization of contact manifolds is known in differential geometry, see \cite{arnold, libermann}. Within a thermodynamics context its use was first advocated in \cite{balian}, and followed up in \cite{entropy}.

\subsection{Homogeneous Hamiltonian dynamics}
Consider as above the cotangent bundle $T^*\mR^3$ {\it minus} its zero section, with coordinates $E,S,V,p_E,p_S,p_V$. For any function $K: T^*\mR^3 \to \mR$, consider the standard Hamiltonian differential equations
\bq
\label{homham}
\begin{array}{rcl}
\dot{E} &=& \frac{\partial K}{\partial p_E}(E,S,V,p_E,p_S,p_V) \\[2mm]
\dot{S} &=& \frac{\partial K}{\partial p_S}(E,S,V,p_E,p_S,p_V) \\[2mm]
\dot{V} &=& \frac{\partial K}{\partial p_V}(E,S,V,p_E,p_S,p_V) \\[2mm]
\dot{p}_E &=& - \frac{\partial K}{\partial E}(E,S,V,p_E,p_S,p_V) \\[2mm]
\dot{p}_S &=& - \frac{\partial K}{\partial S}(E,S,V,p_E,p_S,p_V) \\[2mm]
\dot{p}_V &=& - \frac{\partial K}{\partial V}(E,S,V,p_E,p_S,p_V)
\end{array}
\eq
Now impose the extra condition that $K$ is homogeneous of degree $1$ in $p$, i.e., 
\bq
K(E,S,V,\lambda p_E, \lambda p_S, \lambda p_V)= \lambda K(E,S,V,p_E,p_S,p_V) \mbox{ for all } \lambda \neq 0
\eq
It turns out that for such $K$ the Hamiltonian differential equations \eqref{homham} {\it project} to a {\it contact vector field} on the contact manifold with coordinates (in the energy representation) $E,S,V,T,-P$, and, conversely, that any contact vector field is the projection of such a Hamiltonian dynamics with homogeneous Hamiltonian $K$. The same holds for the entropy representation $E,S,V,\frac{1}{T},\frac{P}{T}$. Generalities concerning this are discussed in the sidebar "\nameref{sb:homogeneous}". The Hamiltonian differential equations for homogeneous $K$ respect the structure of $T^*\mR^3$ (as captured by its Liouville form). Furthermore, such dynamics leaves invariant the homogeneous Lagrangian submanifold $\cL$ (specifying the state properties of the thermodynamic system), if and only if $K$ is zero on on $\cL$. It follows that any Hamiltonian dynamics generated by a function $K$ that is (1) homogeneous of degree $1$ in $p$ and (2) zero on $\cL$, is a {\it feasible dynamics} for the thermodynamic system. This will be the starting point for the definition of {\it port-thermodynamic systems} in the next section.

Interestingly, while the Hamiltonians in the formulation of, e.g., mechanical systems, represent {\it total energy}, the Hamiltonians $K$ as discussed above are {\it dimensionless} (in the sense of dimensional analysis). Furthermore, the {\it contact Hamiltonian} of its projected dynamics (a contact vector field) has dimension of {\it power} in case of the energy representation (with intensive variables $T,-P$), and has dimension of {\it rate of entropy} in case of the entropy representation (with intensive variables $\frac{1}{T}, \frac{P}{T}$). Together with the earlier observation that the dynamics of a thermodynamic system is captured by the dynamics {\it restricted} to the invariant homogeneous Lagrangian submanifold, this once more emphasizes that the Hamiltonian dynamics \eqref{homham} has a rather different interpretation than the Hamiltonian formulation of mechanical (or other physical) systems.

\section{Port-Thermodynamic Systems}
\label{sec:port-thermodynamic}
As argued in the previous section, see \cite{entropy} for further information, any dynamics of a thermodynamic system should respect the geometric structure of the thermodynamic phase space (or its homogeneous symplectic extension). Furthermore, it should leave invariant the Legendre submanifold of the thermodynamic phase space, or equivalently, the corresponding Lagrangian submanifold of the homogeneous symplectic extension. 
Let us adopt the same generality and notation as in the sidebars "\nameref{sb:contact}" and "\nameref{sb:homogeneous}", where the space of all extensive variables is denoted by $\Z$, the symplectic extension by $\T^*\Z$ with coordinates $(z,p)$, and the thermodynamic phase space by $\mP(T^*\Z)$. E.g., in a simple thermodynamic system, $\Z$ is $\mR^3$ with coordinates $E,S,V$ and $\T^*\Z$ has coordinates $E,S,V,p_E,p_S,p_V$, while the coordinates for $\mP(T^*\Z)$ are $E,S,V,T,-P$ (energy representation) or $E,S,V,\frac{1}{T},\frac{P}{T}$ (entropy representation). 
Because of its simplicity, and because it will allow us to define in a natural way {\it ports}, we will focus on the description on the symplectic extension $\T^*\Z$; see \cite{entropy} for details on the resulting projection to the thermodynamic phase space $\mP(T^*\Z)$.

So let us consider a thermodynamic system with constitutive relations (state properties) specified by a homogeneous Lagrangian submanifold $\cL \subset \T^*\Z$. Respecting the geometric structure of the symplectic extension $\T^*\Z$ means that the dynamics is a Hamiltonian vector field $X_K$ on $\T^*\Z$ with $K$ homogeneous of degree $1$ in the $p$-variables. Indeed, by the sidebar "\nameref{sb:homogeneous}" any such a vector field leaves the Liouville form on $\T^*\Z$ invariant.  Furthermore, $X_K$ leaves $\cL$ invariant if and only if $K$ restricted to the homogeneous Lagrangian submanifold $\cL$ is zero. Finally, we will split $K$ into two parts, i.e.,
\bq
K^a + K^cu, \quad u \in \mR^m
\eq
Here $K^a: \T^*\Z \to \mR$ is the homogeneous Hamiltonian corresponding to the autonomous dynamics due to internal non-equilibrium conditions. Next, $K^c=(K^c_1, \cdots, K^c_m)$ is a row vector of homogeneous Hamiltonians (called {\it control} or {\it interaction} Hamiltonians) corresponding to dynamics arising from interaction with the surrounding of the system. This second part of the dynamics is affinely parametrized by a vector $u$ of {\it control}  or {\it input} variables (see however Example \ref{ex:msdint} for a non-affine dependency on $u$). Since $K$ is homogeneous of degree $1$ in $p$ and zero on $\cL$ for all $u\in \mR^m$ this simply means that the $(m+1)$ functions $K^a,K^c_1, \cdots, K^c_m$ are all homogeneous of degree $1$ in $p$ and zero on $\cL$.

By invoking Euler's homogeneous function theorem (cf. Theorem \ref{Euler}) this means that
\bq
\label{KaKc}
\begin{array}{rcll}
K^a &= & p_0\frac{\partial K^a}{\partial p_0} +  p_1\frac{\partial K^a}{\partial p_1} + \cdots + p_n\frac{\partial K^a}{\partial p_n}\\[3mm]
K^c & = & p_0\frac{\partial K^c}{\partial p_0} +  p_1\frac{\partial K^c}{\partial p_1} + \cdots + p_n\frac{\partial K^c}{\partial p_n}
\end{array}
\eq
where the functions $\frac{\partial K^a}{\partial p_i}$, as well as the elements of the $m$-dimensional row vectors of partial derivatives $\frac{\partial K^c}{\partial p_i}$, $i=0,1, \cdots, n$, are all homogeneous of degree $0$ in the $p$-variables. (Note that this does {\it not} necessarily mean that these functions are independent of $p$; although of course this is an important special case.)

There are two more constraints on $K^a$, as imposed by the First and Second Law. Since the energy and entropy variables $E,S$ are among the extensive variables $z_0,z_1, \cdots,z_n$, let us take $E=z_0, S=z_1$.
With this convention, in the internal dynamics $X_{K^a}$ we have $\dot{E}=\frac{\partial K^a}{\partial p_0}$. Hence, by the First Law the energy of the system without interaction with the surrounding (i.e., $u=0$) should be conserved, implying that necessarily $\frac{\partial K^a}{\partial p_0}|_{\cL}=0$.
Similarly, $\frac{\partial K^a}{\partial p_1}$ is equal to $\dot{S}$ in the internal dynamics $X_{K^a}$. Hence by the Second Law of thermodynamics necessarily $\frac{\partial K^a}{\partial p_1}|_{\cL} \geq 0$. 

Such constraints do not hold for the control (interaction) Hamiltonians $K^c$. In fact, the control Hamiltonians may be utilized to define natural {\it output} variables conjugated to the inputs $u$. First option is to define the $m$-dimensional row vector
\bq
y_p=\frac{\partial K^c}{\partial p_0},
\eq
with the subscript $p$ in $y_p$ standing for {\it power}.
Then it follows that along the dynamics $X_K$, with $K=K^a + K^cu$,
\bq
\frac{d}{dt} E= y_pu,
\eq
and thus $y_p$ is the vector of {\it power-conjugate} (passive) outputs corresponding to the input vector $u$. We will call the pair $(u, y_p)$ the {\it power port} of the system. Similarly, by defining the $m$-dimensional row vector ($re$ for 'rate of entropy')
\bq
y_{re}=\frac{\partial K^c}{\partial p_1}
\eq
it follows that along the full dynamics $X_K$
\bq
\frac{d}{dt} S \geq y_{re}u,
\eq
Hence $y_{re}$ is the output vector which is conjugate to $u$ in terms of {\it entropy flow}. The pair $(u, y_{re})$ is called the {\it rate of entropy} port of the system. Note that, in principle, we could also define outputs conjugated to $u$ for other extensive variables as well; thus leading to other ports; see for instance the example of chemical reaction networks below, cf. Example \ref{ex:chem1}. 

We summarize this in the following definition of a {\it port-thermodynamic system} \cite{entropy}.
\begin{definition}
Consider the space of extensive variables $\Z$. A port-thermodynamic system is a pair $(\cL,K)$, where $\cL \subset \T^*\Z$ is a homogeneous Lagrangian submanifold describing the state properties, and $K=K^a + K^cu, u \in \mR^m,$ is a Hamiltonian on $\T^*\Z$, homogeneous of degree $1$ in $p$, and zero restricted to $\cL$.
Let $z=(z_0,z_1, \cdots,z_n)$ with $z_0=E$ (energy), and $z_1=S$ (entropy). Then the power conjugate output is defined as $y_p=\frac{\partial K^c}{\partial p_0}|_{\cL}$, and the rate of entropy conjugate output as $y_{re}=\frac{\partial K^c}{\partial p_1}|_{\cL}$.
\end{definition}

The following examples, mostly taken from \cite{entropy}, illustrate this definition.
\begin{example}[Heat compartment]\label{ex:heatcomp}
Consider a heat compartment, exchanging heat with its surrounding. Its thermodynamic properties are described by the extensive variables $S$ (entropy) and $E$ (internal energy), with $E$ expressed as a function $E=E(S)$ (energy representation) of $S$. Its state properties are given by the homogeneous Lagrangian submanifold
\bq
\cL= \{(S,E,p_S,p_E) \mid E=E(S), p_S=-p_EE'(S) \},
\eq
corresponding to the generating function $-p_EE(S)$. Since there is no internal dynamics, $K^a$ is absent. Hence, taking $u$ as the {\it rate of entropy flow} corresponds to considering the homogeneous Hamiltonian
\bq
K^c=p_S + p_E E'(S), 
\eq
which is clearly zero on $\cL$. This yields on $\cL$ the dynamics, entailing entropy and energy balance
\bq
\begin{array}{rcllcl}
\dot{S} &=& u  \quad & \dot{p}_S &=& -p_EE''(S)u \\[2mm]
\dot{E} &=&E'(S)u \quad & \dot{p}_E &=& 0,
\end{array}
\eq
with power conjugate output $y_p$ equal to the temperature $E'(S)$. 
Defining the homogeneous coordinate $\gamma=-\frac{p_S}{p_E}$ leads to the projected dynamics on the thermodynamic phase space $\mP(T^*\mR^2)$
\bq
\label{heatc}
\begin{array}{rcl}
\dot{S} &=& u\\[2mm]
\dot{E} &=&E'(S)u\\[2mm]
\dot{\gamma} & = & - \frac{\dot{p}_S}{p_E} = E''(S)u
\end{array}
\eq
This is a contact vector field with contact Hamiltonian $\widehat{K}^c= E'(S) - \gamma$, which leaves the Legendre submanifold
\bq
L= \{(S,E,\gamma) \in \mP(T^*\mR^2) \mid E=E(S), \gamma=E'(S) \}
\eq
invariant.
Alternatively, if we take instead the incoming {\it heat flow} as input $v$, then the Hamiltonian is given by
\bq
K=(p_S \frac{1}{E'(S)} + p_E)v,
\eq
leading to the rate of entropy conjugate output $y_{re}$ given by the reciprocal temperature $y_{re}=\frac{1}{E'(S)}$.
\end{example}

\begin{example}[Mass-spring-damper system]\label{ex:msd}
This is an example that normally would not be considered to be a thermodynamic system. Nevertheless, in view of the dissipation of energy due to the damper, there {\it is} a thermodynamic component, namely the heat irreversibly released by the damper. Also, the subsequent example of a gas-piston-damper system will turn out to be analogous.
Consider a mass-spring-damper system in one-dimensional motion, composed of a mass $m$ with momentum $\pi$, linear spring with stiffness $k$ and extension $w$, and linear damper with damping coefficient $d$. In order to take into account the thermal energy and the entropy production arising from the heat produced by the damper, the variables of the mechanical system are augmented with an entropy variable $S$ and internal energy $U(S)$. (For instance, if the system is isothermal, i.e. in thermodynamic equilibrium with a thermostat at temperature $T_0$, the internal energy is $U(S) = T_0S$.) This leads to the total set of extensive variables $w$, $\pi$, $S$, $E=\frac{1}{2}kw^2 + \frac{\pi^2}{2m} + U(S)$ (total energy). The state properties of the system are described by the Lagrangian submanifold $\cL$ with generating function (in energy representation)
\bq
-p_E\left(\frac{1}{2}kw^2 + \frac{\pi^2}{2m} + U(S)\right) ,
\eq
defining the state properties
\bq
\begin{array}{rcl}
\cL & = & \{(w,\pi,S,E,p_z,p_{\pi},p_S,p_E) \mid E=\frac{1}{2}kw^2 + \frac{\pi^2}{2m} + U(S), 
\\[2mm]
&& p_w=-p_Ekw, p_{\pi}= -p_E \frac{\pi}{m},p_S=-p_EU'(S) \}
\end{array}
\eq
The dynamics is given by the following homogeneous Hamiltonian, zero on $\cL$,
\bq
\label{msdK}
K= p_w\frac{\pi}{m} +p_{\pi}\left(-kw -d\frac{\pi}{m}\right) +p_S \frac{d (\frac{\pi}{m})^2}{U'(S)} + \left(p_{\pi} + p_E \frac{\pi}{m}\right)u,
\eq
where $u$ is an external force. The power-conjugate (passive) output $y_p=\frac{\pi}{m}$ is the velocity of the mass.
\end{example}

\begin{example}[Gas-piston-damper system]
\label{ex:gaspiston}
Consider a gas in a thermally isolated cylinder closed by a piston. Assuming the thermodynamic properties of the system to be covered by the properties of the gas, the system is completely analogous to the previous example, replacing $w$ by volume $V$ and the partial energy $\frac{1}{2}kw^2 + U(S)$ by the internal energy $U(S,V)$ of the gas. The dynamics of the gas-piston-damper system, with piston actuated by a force $u$, is given by the Hamiltonian
\bq
K= p_V\frac{\pi}{m} +p_{\pi}\left(-\frac{\partial U}{\partial V} -d\frac{\pi}{m}\right) +p_S \frac{d (\frac{\pi}{m})^2}{\frac{\partial U}{\partial S}} + \left(p_{\pi} + p_E \frac{\pi}{m}\right)u,
\eq
where the power-conjugate output $y_p=\frac{\pi}{m}$ is the velocity of the piston.
\end{example}

\begin{example}[Chemical reaction networks \cite{louvain}]
\label{ex:chem1}
Consider a chemical reaction network in {\it entropy representation}, as in subsection "\nameref{subsec:chemical}", with the entropy $S$ represented as a function $S=S(E,x)$ of the vector of chemical concentrations $x$ and energy $E$. Then the homogeneous Lagrangian submanifold describing the state properties of the reaction network is given as
\bq
\cL = \{(x,S,E,p_x,p_S,p_E) \mid 
S=S(E, x),p_x=-p_S\frac{\partial S}{\partial x}(E, x), p_E=-p_S\frac{\partial S}{\partial E}(E, x) \}
\eq
with $\frac{\partial S}{\partial x}(E,x)= - \frac{\mu}{T}, \frac{\partial S}{\partial E}(E,x)=\frac{1}{T}$. The internal dynamics of the chemical reaction network is generated by the Hamiltonian, homogeneous of degree $1$ in $(p_E,p_S,p_x)$ and zero on $\cL$,
\bq
K^a =  - p_x^\top Z \ML \Exp \, \frac{-Z^\top }{R}\frac{\partial S}{\partial x}(E, x) \, - 
p_S \frac{\partial S}{\partial x^\top }(E, x)Z L \Exp \, \frac{-Z^\top }{R}\frac{\partial S}{\partial x}(E, x)
\eq
Furthermore, the control Hamiltonian
\bq
K^c = p_S\frac{\partial S}{\partial E}(E,x) + p_E ,
\eq
corresponds to a {\it heat flow} input, and a rate of entropy conjugate output $y_e=\frac{\partial S}{\partial E}(x,E)_{|{\mathcal{L}}}$ equal to the reciprocal temperature. 
Another possible choice is
\bq
K^c = p_S\frac{\partial S}{\partial x_i}(E,x) + p_{xi} ,
\eq
corresponding to material in/outflow of the $i$-th chemical species, with rate of entropy conjugate output $y_{re}=\frac{\partial S}{\partial x_i}(E,x)_{|{\mathcal{L}}}$ equal to the chemical potential $\mu_i$ of the $i$-th chemical species divided by the temperature $T$. 
\end{example}
In this last example the internal dynamics is fully {\it irreversible}: the system follows a first-order dynamics and converges to a state where the chemical potentials are equal (very much like {\it consensus dynamics}). In the two examples given before (mass-spring-damper system and gas-piston-damper system) this is of course different. Indeed, although there is irreversible increase of entropy due to the damper action there is also an internal dynamics corresponding to the oscillatory transformation of kinetic energy into potential energy and conversely (second-order dynamics). Thus we conclude that internal dynamics of thermodynamic systems does not necessarily correspond to irreversible dynamics.

In {\it composite}, nonhomogeneous, thermodynamic systems, cf. subsection "\nameref{subsec:basic}", there is typically no {\it single} energy or entropy. In this case the constraints on the internal dynamics are different: the {\it sum} of the energies needs to be conserved, and likewise the {\it sum} of the entropies needs to be increasing. A simple example is the following; see \cite{entropy} for further information.

\begin{example}[Heat exchanger]\label{ex:heatexch}
Consider two heat compartments as in Example \ref{ex:heatcomp}, exchanging a heat flow through a conducting wall according to Fourier's law. The three extensive variables are $S_1,S_2$ (entropies of the two compartments) and $E$ (total internal energy). The state properties are described by the homogeneous Lagrangian submanifold
\bq
\cL  = \{(S_1,S_2,E,p_{S_1},p_{S_2},p_E) \mid E=E_1(S_1) +E_2(S_2), 
p_{S_1}=-p_EE_1'(S_1), p_{S_2}=-p_EE_2'(S_2) \},
\eq
corresponding to the generating function $-p_E\left(E_1(S_1) + E_2(S_2)\right)$, with $E_1,E_2$ the internal energies of the two compartments. Denoting the temperatures $T_1=E_1'(S_1),$ $ T_2=E_2'(S_2)$, the internal dynamics of the two-component thermodynamic system corresponding to Fourier's law is given by the Hamiltonian
\bq
\label{ex:heatexchK}
K^a = \lambda (\frac{1}{T_1} - \frac{1}{T_2})(p_{S_1}T_2 - p_{S_2}T_1),
\eq
with $\lambda$ Fourier's conduction coefficient. The total entropy on $\cL$ satisfies
\bq
\frac{d}{dt}({S}_1 + {S}_2) = \lambda (\frac{1}{T_1} - \frac{1}{T_2})(T_2 - T_1) \geq 0
\eq
\end{example}

\subsection{Equivalent parametrizations of the dynamics}
\label{subsec:equivalent}
As already discussed in the subsection "\nameref{subsec:paradigm}", Gibbs' fundamental relation leads to the consideration of {\it non-minimal} state space representations, involving all the extensive and intensive variables. The minimal state space is the Legendre submanifold $L$ of the thermodynamic phase space, specifying the state properties of the thermodynamic system at hand. In the same way, the definition of a port-thermodynamic system entails dynamics on the {\it whole} thermodynamic phase space (or, equivalently, its symplectic extension), leaving invariant the Legendre submanifold $L$ (or the homogeneous Lagrangian submanifold $\cL$). This means that the dynamics on $L$ can be parametrized in different ways, either by extensive or by intensive variables (or mixtures of them); similar to the use of the various thermodynamical potentials to describe $L$ (see section "\nameref{sec:gibbs}"). Conversely, these different parametrizations of the dynamics are overarched by the dynamics on the whole thermodynamic phase space or its symplectic extension. 

More specifically, consider a port-thermodynamic system with homogeneous Hamiltonian dynamics $X_K$ on $\T^*\Z$ with natural coordinates $(z,p)$, where $K=K^a + K^cu$ is zero on the homogeneous Lagrangian submanifold $\cL$. Let as before, $z_0=E, z_1=S$ and $p_0=p_E,p_1=p_S$. The simplest parametrizations of the dynamics on $\cL$ are obtained by considering the dynamics of the extensive variables $z_1, \cdots,z_n$, corresponding to the energy representation $E=E(z_1, z_2\cdots,z_n)$, or by considering the dynamics of $z_0,z_2, \cdots,z_n$, corresponding to the entropy representation $S=S(z_0,z_2, \cdot, z_n)$. 
On the other hand, the dynamics can equally well be parameterized by considering the dynamics of {\it intensive} variables $\gamma_1, \cdots, \gamma_n$ obtained from the co-extensive variables $p=(p_0,p_1, \cdots,p_n)$. For example, in the energy representation the intensive variables are given as (since $p_0=p_E$)
\bq
\gamma_1= \frac{p_1}{-p_0}, \cdots, \gamma_n= \frac{p_n}{-p_0},
\eq
which restricted to $\cL$ are given by
\bq
\label{leg}
\frac{\partial E}{\partial z_1}(z_1, \cdots,z_n), \cdots, \frac{\partial E}{\partial z_n}(z_1, \cdots,z_n)
\eq
Denote $\widetilde{z}:=(z_1, \cdots,z_n)^\top $, and $\gamma= (\gamma_1, \cdots, \gamma_n)^\top$. Then \eqref{leg} defines a mapping
\bq
\widetilde{z} \mapsto \gamma = \frac{\partial E}{\partial \widetilde{z}}(\widetilde{z})
\eq
Assuming the $n \times n$ Hessian matrix $\frac{\partial^2 E}{\partial \widetilde{z}^2}(\widetilde{z})$ to be {\it invertible} we can then define the {\it Legendre transform} $E^*(\gamma)$ of the function $E(\widetilde{z})$, satisfying the equalities
\bq
\label{leg1}
\gamma= \frac{\partial E}{\partial \widetilde{z}}(\widetilde{z}), \quad \widetilde{z}= \frac{\partial E^*}{\partial \gamma}(\gamma), \quad \frac{\partial^2 E}{\partial \widetilde{z}^2}(\widetilde{z})= \big(\frac{\partial^2 E^*}{\partial \gamma^2}(\gamma)\big)^{-1}
\eq
These equalities allow us to rewrite the dynamics of the $n$ extensive variables $\widetilde{z}$ into dynamics of the $n$ intensive variables $\gamma$, since
\bq
\frac{\partial^2 E^*}{\partial \gamma^2}(\gamma) \dot{\gamma} = \dot{\widetilde{z}}
\eq
A very simple example of this was already provided in Example \ref{ex:heatcomp} (the heat compartment); see \eqref{heatc}. The dynamics on $\cL$ is described by $\dot{S}=u$, or by $\dot{\gamma}=E''(S)u$, where the extensive variable $S$ (entropy) is related to the intensive variable $\gamma $ (temperature) by $\gamma=E'(S)$. A more involved case is Example \ref{ex:gaspiston} (the gas-piston-damper system). Take for simplicity $d=0$ (no damping). Then the dynamics on $\cL$ is given in extensive variables $S,V,\pi$ as
\bq
\begin{array}{rcl}
\dot{S} & = & 0 \\[2mm]
\dot{V} & = & \frac{\pi}{m} \\[2mm]
\dot{\pi} & = & - \frac{\partial U}{\partial V} + u
\end{array}
\eq
On the other hand, the dynamics in the intensive variables $\gamma_1=\frac{\partial U}{\partial S}$ (temperature), $\gamma_2=\frac{\partial U}{\partial V}$ (minus the pressure), and $\gamma_3=\frac{\pi}{m}$ (velocity), can be computed as
\bq
\begin{array}{rcl}
\dot{\gamma}_1 & = & \gamma_3 \frac{\partial^2 U}{\partial S \partial V} \\[2mm]
\dot{\gamma}_2 & = & \gamma_3 \frac{\partial^2 U}{\partial V^2}\\[2mm]
\dot{\gamma}_3 & = & - \frac{1}{m}\gamma_2 + u
\end{array}
\eq
which can be written fully in terms of $\gamma$ by computing the Legendre transform $U^*(\gamma_1,\gamma_2)$ of $U(S,V)$, and using the fact that the Hessian matrix of $U$ is the inverse of the Hessian matrix of $U^*$, cf. \eqref{leg1}. Similar computations can be done in order to obtain a parametrization of the dynamics on $\cL$ in terms of the intensive variables corresponding to the entropy representation. In general, the transformation of the dynamics in extensive variables into the description of the dynamics in intensive variables is similar to the transformation of port-Hamiltonian dynamics in energy variables to its description in {\it co-energy} variables, see e.g. \cite{schjeltsema}. This is also closely related to the (generalized) Brayton-Moser formulation of physical systems \cite{jeltsema,schjeltsema}.

\section{Ports and Interconnections}
The definition of ports enables the interconnection of thermodynamic systems, so as to obtain complex systems from simpler building blocks. This will be only discussed through two examples, and the reader is referred to \cite{entropy} for a more elaborated treatment.
Let us start with the case of {\it power-port} interconnections of port-thermodynamic systems, corresponding to power flow exchange. This is the standard situation in physical network modeling of interconnected systems, in particular in port-based modeling theory, see e.g. \cite{golo}. Consider for simplicity two port-thermodynamic systems, with input vectors $u_1$, respectively $u_2$, and the {\it power-conjugate} outputs $y_{p1},y_{p2}$ as introduced in the definition of port-thermodynamic systems. Then consider interconnection constraints satisfying the {\it power-conservation} property
\bq
y_{p1}^\top u_1 + y_{p2}^\top u_2=0,
\eq
in accordance with the First Law. More generally, in case of an additional {\it external} power port with variables $u,y_p$, consider power-conserving interconnection constraints satisfying
\bq
y_{p1}^\top u_1 + y_{2p}^\top u_2 +y_p^\top u=0
\eq

\begin{example}[Mass-spring-damper system \cite{entropy}]
\label{ex:msdint}
Let us show that the thermodynamic formulation of the system in Example \ref{ex:msd} also results from the interconnection of its three subsystems: mass, spring and damper. The same analysis also applies to the gas-piston-damper system of Example \ref{ex:gaspiston}.\\
I. {\it Mass subsystem} (leaving out irrelevant entropy). The state properties are given by the homogeneous Lagrangian submanifold
\bq
\cL_m= \{(\pi,\kappa,p_{\pi},p_{\kappa}) \mid \kappa=\frac{\pi^2}{2m}, \, p_{\pi} = -p_{\kappa} \frac{\pi}{m} \},
\eq
with energy $\kappa$ (kinetic energy), and dynamics generated by the Hamiltonian
\bq
K_m= (p_{\kappa} \frac{\pi}{m} +p_{\pi})u_m,
\eq
corresponding to $\dot{\pi}=u_m, y_m=\frac{\pi}{m}$.\\
II. {\it Spring subsystem} (again leaving out irrelevant entropy). The state properties are given by
\bq
\cL_s= \{(w,P,p_w,p_P) \mid P=\frac{1}{2}kq^2, \,p_w = -p_P kw \},
\eq
with energy $P$ (spring potential energy), and dynamics generated by the Hamiltonian
\bq
K_s= (p_P kw +p_{w})u_s,
\eq
corresponding to $\dot{w}=u_s, y_s=kz$.\\
III. {\it Damper subsystem}.
The state properties are given by
\bq
\cL_d=\{ (S,U) \mid U=U(S), \, p_S=-p_UU'(S) \},
\eq
involving the entropy $S$ and an internal energy $U(S)$, while the dynamics is generated by the Hamiltonian
\bq
K_d= (p_U + p_{S} \frac{1}{U'(S)})du_d^2
\eq
with $d$ the damping constant, and power-conjugate output $y_d:= du_d$ equal to the damping force.\\
Now interconnect the three subsystems via their power-ports $(u_m,y_m), (u_s,y_s), (u_d,y_d)$:
\bq
u_m= -y_s - y_d, \; u_s=y_m=u_d
\eq
This results (after setting $p_{\kappa}=p_P=p_U$) in the interconnected port-thermodynamic system with total Hamiltonian $K_m + K_s + K_d$ 
given by the Hamiltonian for $u=0$ as obtained before in Example \ref{ex:msd}, eqn. \eqref{msdK}.
\end{example}

The situation for interconnection via rate-of-entropy ports is {\it different}, since by the Second Law the rate of total entropy is not necessarily zero, but greater or equal than zero. This is illustrated by the following example.
\begin{example}
\label{ex:heatexch1}
The heat exchanger as in Example \ref{ex:heatexch} can be modelled as the interconnection of two heat compartments as in Example \ref{ex:heatcomp}, via their rate of entropy ports $(v_i, y_{rei})$, where $y_{rei}=\frac{1}{E'(S_i)}$, $i=1,2$. The interconnection is defined as
\bq
v_1=-v_2= \lambda (\frac{1}{y_{re2}} - \frac{1}{y_{re1}}),
\eq
with $\lambda >0$ Fourier's conduction coefficient. Clearly, this interconnection, corresponding to a conducting wall, is not entropy conserving, but instead corresponds to increase of entropy
\bq
y_{re1}v_1 +y_{re2}v_2= \lambda \big(\frac{1}{y_{re2}} - \frac{1}{y_{re1}} \big)(y_{re1}-y_{re2} ) \geq 0 
\eq
\end{example}
%

%

\section{Conclusions}
Emphasis in this paper has been put on two aspects: a clear {\it cyclo-dissipativity} interpretation of classical thermodynamics, and, starting from this, a {\it geometric} (coordinate-free) formulation of the state properties of a thermodynamic system and its dynamics through contact and symplectic geometry. Both aspects are considered to be essential in aligning thermodynamics with modern dynamical systems and control theory, and for integration of thermodynamics in unified frameworks for complex systems modeling and control.

Many other aspects of thermodynamics have not been covered at all. Among them, the maximum entropy principle (using the concavity of the entropy function), the stability analysis of forced equilibria of nonequilibrium thermodynamic systems by minimal irreversible entropy production \cite{kondepudi}, and mesoscopic thermodynamics, see e.g. \cite{grmela}, seem to be especially relevant for systems and control. 
Finally, as mentioned before, a major challenge lies in the connection of thermodynamics with information theory, aimed at uniting control strategies based on energy shaping and energy routing, with an information processing point of view. 

\subsection{Acknowledgements} I thank Bernhard Maschke for an inspiring and enduring collaboration; without this the writing of this paper would have been unthinkable.

\section{Sidebar: Cyclo-Dissipativity Theory}
\label{sb:cyclo}

Dissipativity theory originates from the seminal work of Willems \cite{willems1972}, continued by Hill \& Moylan (e.g. \cite{hillmoylan1980}) and others; see \cite{passivitybook} for an updated and extended exposition. The notion of dissipativity was relaxed to {\it cyclo-dissipativity} in \cite{willems1973}, and further explored in the report \cite{hillmoylan1975}, with recent extensions in \cite{vdscyclo}. In the present context of thermodynamics it is appropriate to emphasize {\it cyclo}-dissipativity. 

Consider a system with vector of state variables $x \in \X$, and vector of external variables $w \in \W$ (comprising the inputs $u$ and outputs $y$). Furthermore, consider a {\it supply rate} $s:\W \to \mR$.
\begin{definition}\label{cyclo-ext}
A system is {\it cyclo-dissipative} (for supply rate $s$) if 
\begin{equation}\label{cyclic}
\int\limits_{t_1}^{t_2}  s\big(w(t)\big) dt \geq 0
\end{equation}
for all $t_2\geq t_1$ and all external trajectories $w(\cdot)$ such that $x(t_2)=x(t_1)$. In case \eqref{cyclic} holds with equality, the system is {\it cyclo-lossless}. 
Furthermore, the system is {\it cyclo-dissipative with respect to }$x^*$ if \eqref{cyclic} holds for all $t_2\geq t_1$ and all external trajectories $w(\cdot)$ such that $x(t_2)=x(t_1)=x^*$. Finally, it is {\it cyclo-lossless with respect to} $x^*$ if this holds with equality.
\end{definition}
Interpreting $s(w)$ as 'power' provided by the surrounding to the system, cyclo-dissipativity means that for any cyclic trajectory the net amount of 'energy' supplied to the system is non-negative, and zero in case of cyclo-losslessness. Hence a cyclo-dissipative system cannot generate (but only dissipate) 'energy'. In the special case of the supply rate $s(u,y)=y^Tu$, with the equally dimensioned vectors of inputs $u$ and outputs $y$ comprising the vector of external variables $w$, 'cyclo-dissipativity' is referred to as 'cyclo-passivity'.

Apart from the requirement that the state of the system follows a cyclic process, Definition \ref{cyclo-ext} entails an {\it external} characterization of cyclo-dissipativity. In order to relate this external characterization to the internal state dynamics the notions of {\it dissipation (in)equality} and {\it storage function} are introduced.

\begin{definition}
Consider a system with state vector $x$, vector of external variables $w$, and supply rate $s$. 
A function $F:\X \to \mR$ is called a storage function if it satisfies the dissipation inequality
\begin{equation}\label{diss}
F\big(x(t_2)\big) - F\big(x(t_1)\big) \leq \int\limits_{t_1}^{t_2} s(w(t)) dt 
\end{equation}
for all $t_2\geq t_1$, all initial conditions $x(t_1)$, and all external trajectories $w(\cdot)$, where $x(t_2)$ is the state at time $t_2$ corresponding to initial condition $x(t_1)$ and external trajectory $w(\cdot)$.
Eq. \eqref{diss} with equality is called the {\it dissipation equality}.
\end{definition}

Interpreting as before $s\big(w(t)\big)$ as 'power' supplied to the system at time $t$, and $F\big(x(t)\big)$ as stored 'energy' while the system is at state $x(t)$, the existence of a storage function means that increase of the stored energy can only occur due to externally supplied power. The following theorem \cite{vdscyclo} extends the results in \cite{hillmoylan1975} and shows the equivalence between, on the one hand, the external characterization of cyclo-dissipativity and cyclo-losslessness, and, on the other hand, the existence of storage functions. 

\begin{theorem}
\label{th:cyclo}
Consider a system with supply rate $s$. If there exists a storage function $F$ then the system is cyclo-dissipative, and it is cyclo-lossless if $F$ satisfies \eqref{diss} with equality. Assume that the system is reachable from some 
ground state $x^*$ and controllable to this same state $x^*$. [It is immediate that this property is independent of the choice of $x^*$.] Define the (possibly extended) functions $F_{ac}: \X \to \mR \cup \infty$ and $F_{rc} : \X \to - \infty \cup \mR$ as
\begin{equation}
\begin{aligned}
F_{ac}(x) &= \!\!\!\! \mathop{\sup_{w,\T\geq 0 \,|}}_{x(0)=x,\, x(\T)=x^*} - \int\limits_0^\T s(w(t)) dt, \\[2mm]
F_{rc}(x) &= \!\!\!\!\!\! \mathop{\inf_{w,\T\geq 0 \,|}}_{x(-\T)=x^*,\, x(0)=x}   \ \ \int\limits_{-\T}^0s(w(t)) dt,
\end{aligned}
\end{equation}
where the supremum and infimum are taken over all external trajectories $w(\cdot)$ and $\T\geq 0$, satisfying $x(0)=x,\, x(\T)=x^*$, respectively $x(-\T)=x^*,\, x(0)=x$.
Then the system is cyclo-dissipative with respect to $x^*$ if and only if
\begin{equation}\label{acrc}
F_{ac}(x) \leq F_{rc}(x), \ \text{for all} \ x \in \X.
\end{equation}
Furthermore, if the system is cyclo-dissipative with respect to $x^*$, then $F_{ac}: \X \to \mR$ and $F_{rc}: \X \to \mR$, and they define storage functions, implying the system is cyclo-dissipative. Furthermore, $F_{ac}(x^*) = F_{rc}(x^*)=0$, while any other storage function $F$ satisfies
\begin{equation*}
F_{ac}(x) \leq F(x) - F(x^*) \leq F_{rc}(x)
\end{equation*}
If the system is cyclo-lossless with respect to $x^*$ then $F_{ac}(x)= F_{rc}(x), x \in \X$, implying uniqueness (up to a constant) of the storage function.
\end{theorem}

Note that the first statement of this theorem, existence of a storage function implying cyclo-dissipativity, is obvious: simply substitute $x(t_1)=x(t_2)$ into \eqref{diss}. For the proof of the rest of this theorem we refer to \cite{vdscyclo}.

In general, apart from the cyclo-lossless case, storage functions of dissipative systems are far from unique. The following proposition from \cite{vdscyclo} ensures uniqueness by imposing a {\it weakened} form of cyclo-losslessness.
\begin{proposition}
\label{prop:revers}
Suppose the system is reachable from and controllable to $x^*$ and cyclo-dissipative with respect to $x^*$. Assume additionally that for every $x$ there exists a solution $(x_l(\cdot), w_l(\cdot))$ on some time-interval $[0,\T_l]$ such that $x_l(0)=x_l(\T_l)=x^*$ and $x_l(\tau)=x$ for some $\tau \in [0,\T_l]$, satisfying
\bq
\label{weakl}
\int_{0}^{\T_l} s(w_l(t)) dt =0
\eq
Then $F_{ac}(x)=F_{rc}(x)$ for all $x \in \X$, and the storage function is unique up to a constant, and given by
\bq
F(x)=\int_{0}^{\tau} s(w_l(t)) dt = - \int_{\tau}^{\T_l} s(w_l(t)) dt,
\eq
where $x_l(\tau)=x$ and $x_l(0)=x_l(\T_l)=x^*$.
\end{proposition}

Note that \eqref{weakl} means that the system is weakly cyclo-lossless with respect to $x^*$, in the sense that for every $x$ there exists {\it at least one} cyclic trajectory passing through $x$ and $x^*$ satisfying \eqref{cyclic} with {\it equality}; while other cyclic trajectories passing through $x$ and $x^*$ satisfy \eqref{cyclic} only with {\it inequality}.

The stronger notion of dissipativity \cite{willems1972}, originally introduced {\it before} the notion of cyclo-dissipativity, historically {\it starts} from the dissipation inequality \eqref{diss}, restricting to {\it non-negative} storage functions. See also \cite{willems2007} for further information.

\begin{definition}
Consider a system with state vector $x$, vector of external variables $w$, and supply rate $s$.
The system is {\it dissipative} (for the supply rate $s$) if there exists a {\it non-negative} storage function $F$. [Since addition of an arbitrary constant to a storage function again leads to a storage function, the requirement of non-negativity of $F$ can be relaxed to $F$ being {\it bounded from below}.] Furthermore, it is {\it lossless} if there exists a non-negative storage function $F$ satisfying the dissipation inequality \eqref{diss} with {\it equality}.
\end{definition}
Non-negative storage functions are candidate Lyapunov functions for the internal state dynamics for $s(w)=0$; see \cite{willems1972, hillmoylan1980, passivitybook}. In this way, dissipativity theory connects the external stability properties of the system to internal, Lyapunov, stability theory.
An {\it external} characterization of dissipativity is obtained as follows \cite{willems1972}. 
%
\begin{theorem}
\label{tm:diss}
The system is dissipative (for the supply rate $s$) if and only if 
\begin{equation}
F_a(x):= \sup_{w, \T\geq 0} - \int\limits_0^\T s\big(w(t)\big) dt < \infty
\end{equation}
for every $x$, where the supremum is taken over all external trajectories $w(\cdot)$ of the system corresponding to initial condition $x(0)=x$, and all $\T\geq 0$. Obviously, $F_a(x)\geq 0$. 
Furthermore, if $F_a(x) < \infty$ for every $x$ then $F_a$ is a non-negative storage function, and is in fact the {\it minimal} non-negative storage function. 
If additionally the system is reachable from some ground state $x^*$ then it is dissipative if and only if $F_a(x^*)< \infty$.
\end{theorem}

Interpreting again $s(w)$ as the 'power' supplied to the system, $F_a(x)$ is the {\it maximal 'energy'} that can be extracted from the system at initial condition $x$. Thus Theorem \ref{tm:diss} states that the system is dissipative if and only if from any initial state $x$ only a {\it finite} amount of 'energy' can be extracted. This should be contrasted with the external characterization \eqref{cyclic} of cyclo-dissipativity, stating that the system is cyclo-dissipative if and only if the net 'energy' supplied to the system along any cyclic trajectory is $\geq 0$. In fact, in case the system is cyclo-dissipative system it may still be possible to extract an {\it infinite} amount of 'energy' (namely, if the storage function is not bounded from below). Storage functions that are not bounded from below are not uncommon in physical systems modeling. For example, the gravitational energy between two masses is proportional to $-\frac{1}{r}$, with $r\geq 0$ the distance between the two masses, and is thus {\it not} bounded from below. 

\begin{remark}
For reachable linear systems with quadratic supply rates, dissipativity is often equivalent to $F_a(0)=0$, yielding the familiar external characterization $\int\limits_0^\T s\big(w(t)\big) dt \geq 0$ for all trajectories starting from $x(0)=0$; see e.g. \cite{passivitybook}.
\end{remark}
Finally, by assuming {\it differentiability} of the storage function, the dissipation inequality can be replaced by an (easier) differential version. Consider e.g. an input-state-output system 
\begin{equation}
\begin{aligned}
\dot{x} & = f(x,u),\\
y & = h(x,u),
\end{aligned}
\end{equation}
with state $x \in \X$, and vector of external variables $w=(u,y)$, where $u$ is the vector of inputs, and $y$ the vector of outputs. Consider a supply rate $s(u,y)$. Then a differentiable function $F:\X \to \mR$ satisfies the dissipation inequality \eqref{diss} (and thus is a storage function) if and only if it satisfies the {\it differential dissipation inequality} \cite{willems1972, hillmoylan1980, passivitybook}
\begin{equation}
\label{diss-eq}
\frac{\partial F}{\partial x}(x)f(x,u) \leq s\big(u,h(x,u)\big), \ \text{ for all } x,u,
\end{equation}
while it satisfies the dissipation {\it equality} if \eqref{diss-eq} holds with equality.

\section{Sidebar: Carnot and Caloric Theory}
\label{sb:caloric}
Interestingly, it seems \cite{kondepudi,fowler} that Sadi Carnot started his investigations into the maximal efficiency of steam engines based on the {\it caloric theory}. Scientists before him like Benjamin Franklin \cite{fowler} believed that heat was flowing through material by some (almost) weightless caloric fluid, and that the amount of caloric fluid was {\it conserved}. In the same spirit Carnot's initial idea \cite{fowler} was that just as water flows downhill, caloric fluid flows from hot to cold, and that the steam engine utilizes this caloric flow to produce work, just as a water wheel takes energy from falling water. As a consequence, Carnot originally believed that in his 'Carnot cycle' the amount of heat $Q_h$ absorbed from the hot reservoir is equal to the amount of heat $-Q_c$ released to the cold reservoir. Only later it seems he realized the fallacy of this idea \cite{kondepudi}. By the time his work was made public (only in 1878; although Clapeyron used Carnot's ideas in his description of the Carnot cycle in 1834, after Carnot's death in 1832) the First Law was already fully accepted, notably through the work of Joule, Helmholtz and Mayer \cite{kondepudi}.

From a mathematical point of view (without worrying about the physics), let us suppose that the caloric fluid {\it is} conserved, and that there exists a function $Q$ of the state of the thermodynamic system ($Q$ being the amount of caloric stored in the system) such that $\frac{d}{dt}Q =q$, where $q$ is the caloric (heat) flow.
Then, if additionally the First Law $\frac{d}{dt}E =q -Pu_V$ holds, this would imply the existence of a third function $W$ of the state, defined as $W:= E-Q$, satisfying $\frac{d}{dt}W =- Pu_V$. Said otherwise, the energy $E$ would be the sum of two functions $Q$ and $W$, which are storage functions for the supply rates $q$, respectively $- Pu_V$. Typically, this can only be the case if the thermodynamic system consists of two {\it separate} parts, one for heat storage and one purely mechanical.

\section{Sidebar: Other Views on Entropy}
\label{otherentropy}

While the section "\nameref{sec:secondlaw}" discusses the classical way of defining entropy, going back to Clausius and reinforced by cyclo-dissipativity theory, there are alternative approaches as well. One of them was initiated by Carath\'eodory, and advocated by Born. The basic idea is as follows; see \cite{bamberg} for further information. Consider a simple thermodynamic system (the argument can be quite easily extended to more complicated situations). By the First Law we know that there exists a function $E(x)$ of the state $x$ of the thermodynamic system, representing the stored energy. Now let us consider the one-form
\bq
\beta:= dE +PdV
\eq
We know from the considerations in the sidebar "\nameref{sb:caloric}" that in general $\beta$ is {\it not} an exact one-form; i.e., there does not exist a state function $Q$ such that $\beta=dQ$. However, we can proceed as follows. 
Curves on the state space whose tangent (velocity) vectors at every point of the curve are in the kernel of $\beta$ will be called {\it adiabatic curves}. (Note that by the First Law $\beta$ evaluated at a tangent vector indeed equals the heat flow $q$; justifying the terminology 'adiabatic'.)
Now replace the formulation of the Second Law of thermodynamics as given by Kelvin (see the section "\nameref{sec:secondlaw}") by the following {\it alternative} statement:

{\it Near any state $x$ there exist arbitrarily close states which cannot be joined to $x$ by an adiabatic curve.}

Then by Carath\'eodory's theorem on one-forms, see e.g. \cite{bamberg}, it follows that, although $\beta$ is not exact, there exist functions of the state, called $S$ (entropy) and $\tau$ such that $\beta=\tau dS$. Equating $\tau$ with the absolute temperature $T$ then yields Gibbs' fundamental relation, expressed by saying that the one-form
\bq
dE + PdV - TdS
\eq
is zero on the state space. Note however the differences with the definition of entropy by Clausius, as exposed in "\nameref{sec:secondlaw}": (1) the alternative formulation of the Second Law is different from the classical formulation of the Second Law as expressed by Kelvin (or the equivalent formulation given by Clausius himself \cite{fermi}), (2) heat flow in the Born-Carath\'eodory approach to thermodynamics is a {\it derived} concept; in contrast with the theory of Clausius and cyclo-dissipativity theory, (3) irreversible thermodynamics is not covered.

The most fundamental definition of entropy, but outside the realm of classical, macroscopic, thermodynamics, is the one given in statistical thermodynamics by Ludwig Boltzmann. This definition can be motivated, in a very rudimentary way, as follows. Obviously, $\frac{d}{dt}S \geq \frac{q}{T}$ implies that the entropy of a thermodynamic system without external heat flow can only increase. On the other hand, in accordance with statistical considerations, from a {\it microscopic} point of view it is plausible that the (very high-dimensional) state of an isolated system will converge to the state of {\it highest probability}. This led Boltzmann to establish his fundamental relationship $S=k \log \pi$, where $\pi$ denotes the number of microscopic states corresponding to the macroscopic thermodynamic state, and $k$ is called the Boltzmann constant. In fact, see e.g. \cite{fermi}, in order to demonstrate this fundamental relationship the essential step is to show that $S=f(\pi)$ for a certain function $f$. Namely, consider a system composed of two parts, with entropies $S_1$ and $S_2$, and $\pi_1,\pi_2$ numbers of microscopic states corresponding to the macroscopic states of both parts. Then the entropy of the total system is given by $S=S_1+S_2$, while the number of microscopic states corresponding to the macroscopic thermodynamic state of the total system is $\pi_1 \pi_2$. This implies that the function $f$ should satisfy $f(\pi_1 \pi_2)= f(\pi_1) + f(\pi_2)$, which is rather easily seen to imply that $f(\pi)=k \log \pi$ modulo a constant.
Note that Boltzmann's definition of entropy is derived from an {\it isolated} system point of view. In contrast, the definition of entropy made by Clausius, which can be naturally interpreted from a cyclo-dissipativity point of view, is crucially based on the {\it interaction} of the macroscopic thermodynamical system with heat sources. This dichotomy provokes many stimulating questions. Importantly, Boltzmann's definition of entropy inspired the definition of entropy in information theory as given by Shannon. Control often has two complementary aspects: energy storage and routing of, and information gathering and processing. It is tempting to assume that thermodynamics may provide the key to unify both aspects. 

Another view on entropy in macroscopic thermodynamics, influenced by Boltzmann's statistical dynamics definition, was advocated by Callen \cite{callen}, and followed up by many others. The Second Law of thermodynamics is replaced by the following {\it postulate} about the existence of the entropy as a function of the state and the {\it entropy maximum principle}:

{\it There exists a function (called the entropy) of the extensive variables of any composite system, defined for all states and having the following property: The values assumed by the extensive variables in the absence of an internal constraint are those that maximize the entropy over the manifold of constrained states.}

A basic illustration of this postulate is a composite system consisting of two parts, with energies $E_1$ and $E_2$, constrained by the requirement that $E_1 + E_2$ is constant, together with an internal constraint that the wall between the two parts is non-conducting. Maximization of $S$ over all those $E_1$ and $E_2$ yields the values of $E_1$ and $E_2$ that obtain when the two parts are connected by a conducting wall.

Finally \cite{haddad}, see also \cite{willemsreview} for a review, presents a {\it middle-ground} theory of thermodynamics; i.e., a foundational framework in between macroscopic and statistical thermodynamics. It is based on deterministic large-scale dynamical systems theory and dissipativity theory. In particular, it is aimed at defining entropy and related notions in a rigorous way, making use of equipartition concepts.

\section{Sidebar: Contact Geometry}
\label{sb:contact}

As discussed in section "\nameref{sec:gibbs}", the state properties of a simple thermodynamic system with extensive variables $E,S,V$ and intensive variables $T, P$ are described  by a $2$-dimensional submanifold $L \subset  \mR^5$, which is such that the Gibbs form
\bq
dE - TdS +PdV
\eq
is zero restricted to $L$ (i.e., at any point of $L$ the Gibbs form annihilates every tangent vector to $L$ at this point). This is an example of {\it contact geometry}.
In general \cite{arnold, libermann}, a {\it contact manifold} $M$ is an odd-dimensional manifold endowed with a {\it contact form} $\theta$. Without going into any detail, a one-form $\theta$ on a $(2n+1)$-dimensional manifold $M$ is a contact form if and only around any point in $M$ we can find coordinates $(z_0,z_1, \cdots,z_n, \gamma_1, \cdots, \gamma_n)$ for $M$ such that
\bq
\label{darboux}
\theta = dz_0 - \sum_{k=1}^n \gamma_k dz_k
\eq
In particular, the Gibbs form $dE - TdS + PdV$ is a contact form on $M=\mR^5$.
\begin{remark} In the actual definition \cite{arnold} of a contact manifold $\theta$ only needs to be defined {\it locally}. What counts is the {\it contact distribution}, the $2n$-dimensional subspace of the tangent space at any point of $M$ defined by the kernel of the contact form $\theta$ at this point.
\end{remark}
A {\it Legendre submanifold} of the contact manifold $(M,\theta)$ is a submanifold of maximal dimension restricted to which the contact form $\theta$ is zero. The dimension of any Legendre submanifold of a $(2n+1)$-dimensional contact manifold is equal to $n$. In particular, the state space manifold of a simple thermodynamic system is a $2$-dimensional Legendre submanifold of $(\mR^5, dE - TdS +PdV)$.
Any Legendre submanifold can be locally represented by a {\it generating function} (and usually in many ways).
\begin{proposition} Consider a contact manifold $(M, \theta)$, with $\theta$ being locally given by \eqref{darboux}. Then there exists a partitioning $\{1, \cdots, n\}= I \cup J$, with $I \cap J = \emptyset$, and locally a generating function $F(z_I,\gamma_J)$, with $z_I$ denoting the coordinates $z_i$ with $i \in I$, and $\gamma_J$ denoting the coordinates $\gamma_j$ with $j \in J$, such that $L$ is given as
\bq
\label{legendre}
L= \{(z_0,z_1, \cdots,z_n, \gamma_1, \cdots, \gamma_n) \mid z_0 = F- \gamma_J \frac{\partial F}{\partial \gamma_J}, 
z_J= - \frac{\partial F}{\partial \gamma_J}, \gamma_I=  \frac{\partial F}{\partial z_I} \} ,
\eq
while conversely any submanifold as in \eqref{legendre} is a Legendre submanifold. Furthermore, the different generating functions can be obtained from each other by {\it partial Legendre transform}.
\end{proposition}
In a thermodynamic system the different possible choices of $F(z_I,\gamma_J)$ exactly correspond to the thermodynamic potentials, as discussed in section "\nameref{sec:gibbs}". For example, in the energy representation of a simple thermodynamic system the generating function is the energy $E(S,V)$, with $J$ void, while $F(T,V)$ is the Helmholtz free energy.
This can be immediately extended to the one-form
\bq
dE -TdS +PdV - \sum_{k=1}^m\mu_kdN_k,
\eq
where $N_k$ is the mole number of the $k$-th chemical species and $\mu_k$ its chemical potential, and to more general one-forms corresponding to other physical cases \cite{kondepudi}.

Finally, a vector field $X$ on a contact manifold $(M, \theta)$ is called a {\it contact vector field} if (with $\mL_X$ denoting {\it Lie derivative} with respect to the vector field $X$)
\bq
\mL_X \theta = \rho \theta ,
\eq
for some scalar function $\rho$ on $M$. (Note that this means that $X$ leaves the {\it contact distribution} invariant.) The function $\theta(X)$ is called the {\it contact Hamiltonian} of the contact vector field. Conversely, to any function on $M$ there corresponds a contact vector field. The expression of a contact vector field in local coordinates is somewhat complicated \cite{bravetti, eberard}. Instead we will focus on the, easier, homogeneous Hamiltonian vector fields on the symplectic extension of the contact manifold, which {\it project} to contact vector fields; see the sidebar "\nameref{sb:homogeneous}".


\section{Sidebar: Homogeneous Symplectic Geometry}
\label{sb:homogeneous}

Here we indicate how {\it contact geometry} (as briefly described in the sidebar "\nameref{sb:contact}") can be formulated as {\it homogeneous symplectic geometry}, by adding one dimension. This correspondence is known in differential geometry, see e.g. \cite{arnold, libermann}. Its relevance for the geometric description of thermodynamics was first advocated in \cite{balian}, and followed up in \cite{entropy}.

We start with the collection of all {\it extensive variables}, that is $E,S$, together with the remaining extensive variables such as $V, N_1, \cdots, N_m, \cdots$. The vector of all extensive variables will be denoted by $z=(z_0,z_1, \cdots,z_n) \in \Z$, with $\Z$ the manifold of extensive variables. Next we consider the cotangent bundle $\T^*\Z$ without its zero-section. Given the coordinates $z$ for $\Z$ there are natural coordinates for the cotangent space denoted by $p=(p_0,p_1, \cdots,p_n)$, leading to natural coordinates $(z,p)=(z_0,\cdots,z_n, p_0, \cdots,p_n)$ for $\T^*\Z$. In the case of thermodynamics
\bq
z=(E,S,V, N_1, \cdots,N_m, \cdots), \quad p=(p_E,p_S,p_V, p_{N_1}, \cdots, p_{N_m}, \cdots)
\eq
$\T^*\Z$ is endowed with a natural one-form $\alpha$ (called the Liouville form), in the above coordinates $(z,p)$ given as
\bq
\label{Liouville}
\alpha= p_0dz_0 + p_1dz_1 + \cdots + p_ndz_n
\eq
Then for each $z\in \Z$ and each cotangent space $T^*_z\Z$ we consider the {\it projective space} $\mP(T^*_z\Z)$, given as the set of rays in $T^*_z\Z$, that is, all the non-zero multiples of a non-zero cotangent vector. Thus the projective space $\mP(T^*_z\Z)$ has dimension $n$, and there is a canonical projection $\pi_z: \T^*_z\Z \to \mP(T^*_z\Z)$, where $\T^*_z\Z$ denotes the cotangent space {\it without} the zero vector. The fiber bundle of the projective spaces $\mP(T^*_z\Z)$, $z \in \Z$, over the base manifold $\Z$ will be denoted by $\mP(T^*\Z)$, and defines a {\it contact manifold} of dimension $2n+1$ (one less than the dimension of $\T^*\Z$); cf. \cite{arnold,entropy}. Informally, whenever $p_0 \neq 0$ we divide the Liouville form $\alpha$ in \eqref{Liouville} by $-p_0$, so as to obtain the contact form
\bq
\theta= dz_0 -\gamma_1dz_1  \cdots -\gamma_n dz_n, \quad \gamma_i:= \frac{p_i}{-p_0}
\eq
Furthermore, if $p_0$ happens to be zero then we divide by another $-p_i \neq 0$; just as in the transition from energy to entropy representation. The contact manifold $\mP(T^*\Z)$ defines the {\it canonical thermodynamic phase space}.

Since each cotangent space (minus the zero vector) $\T^*_z\Z$  projects under $\pi_z$ to the projective space $\mP(T^*_z\Z)$, this defines a total projection $\pi: \T^*\Z \to \mP(T^*\Z)$.
All the relevant objects on the contact manifold $\mP(T^*\Z)$ (such as functions, Legendre submanifolds and contact vector fields) can be shown to correspond to objects on $\T^*\Z$ with an added property of {\it homogeneity} in the co-tangent variables $p$, in such a way that they project under $\pi$ to an object on the contact manifold $\mP(T^*\Z)$. 

Let us start with homogeneity of {\it functions}, characterized by Euler's theorem.
\begin{definition}
Let $r \in \mZ$. A function $K: \T^*\Z \to \mR$ is called homogeneous of degree $r$ (in $p$) if
\bq
K(z, \lambda p) = \lambda^r K(q, p), \quad \mbox{ for all } \lambda \neq 0
\eq
\end{definition}
\begin{theorem}[Euler's homogeneous function theorem]\label{Euler}
A differentiable function $K: \T^*\Z \to \mR$ is homogeneous of degree $r$ (in $p$) if and only if
\bq
\sum_{i=0}^{n} p_i\frac{\partial K}{\partial p_i}(z,p)= rK(z,p), \quad \mbox{ for all } (z,p) \in \T^*\Z
\eq
Furthermore, if $K$ is homogeneous of degree $r$, then its derivatives $\frac{\partial K}{\partial p_i}(z,p),$$ i=0,1,\cdots,n,$ are homogeneous of degree $r-1$.
\end{theorem}
Obviously, a function $K: \T^*\Z \to \mR$ that is homogeneous of degree $0$ in $p$ projects to a function on the thermodynamic phase space $\mP(T^*\Z)$. Next we will consider {\it homogeneous Lagrangian submanifolds}. Recall that a Lagrangian submanifold $\cL \subset \T^*\Z$ is a maximal submanifold of $\T^*\Z$ such that $d \alpha |_{\cL}=0$. A Lagrangian submanifold is called {\it homogeneous} if whenever $(z,p) \in \cL$ then also $(z,\lambda p) \in \cL$ for any $0 \neq \lambda \in \mR$.
\begin{proposition}[\cite{entropy}]
A submanifold $\cL \subset \T^*\Z$ is a homogeneous Lagrangian submanifold if and only if $\alpha |_{\cL}=0$.
For any homogeneous Lagrangian submanifold $\cL \subset \T^*\Z$ there exists a Legendre submanifold $L \subset \mP(T^*\Z)$ such that $\cL=\pi^{-1}(L)$, and conversely for any Legendre submanifold $L \subset \mP(T^*\Z)$ there exists a homogeneous Lagrangian submanifold $\cL \subset \T^*\Z$ such that $\cL=\pi^{-1}(L)$.
\end{proposition}

Thus Legendre submanifolds $L$ of the contact manifold $\mP(T^*\Z)$ (the canonical thermodynamic phase space) correspond to homogeneous Lagrangian submanifolds $\cL$ of $\T^*\Z$. Furthermore, let $F(z_I,\gamma_J)$, with $\{1,\cdots,n\}=I \cup J$, $I \cap J = \emptyset$, be a generating function for the Legendre submanifold $L \subset \mP(T^*\Z)$. Then a generating function for the corresponding homogeneous Lagrangian submanifold $\cL \subset \T^*\Z$ (such that $\cL=\pi^{-1}(L)$) is given by
\bq
G(z_I,p_0,p_J) = -p_0F(z_I,\frac{p_J}{-p_0}),
\eq
in the sense that 
\bq
\cL = \{(z,p) \mid q_0= - \frac{\partial G}{\partial p_0}, \, q_J=- \frac{\partial G}{\partial p_J}, \, p_I= \frac{\partial G}{\partial q_I} \}
\eq
Finally, we come to {\it dynamics}. Recall that for any function $K: \T^*\Z \to \mR$ the \emph{Hamiltonian vector field} $X_K$ on $\T^*\Z$ is defined by the standard Hamiltonian equations
\bq
\dot{z}_i = \frac{\partial K}{\partial p_i}(z,p), \quad \dot{p}_i = - \frac{\partial K}{\partial z_i}(z,p), \quad i=0,1\cdots, n
\eq
Now impose on the Hamiltonians $K: \T^*\Z \to \mR$ the condition that they are homogeneous of degree $1$ in $p$, i.e., $K(q,\lambda p)= \lambda K(q,p)$ for all $\lambda \neq 0$.
\begin{proposition}
If $K: \T^*\Z \to \mR$ is homogeneous of degree $1$ in $p$ then its Hamiltonian vector field $X_K$ is such that ($\mL_X$ denoting Lie derivative with respect to the vector field $X$)
\bq
\mL_{X_K} \alpha =0
\eq
Conversely, if $\mL_{X}\alpha=0$, then $X=X_K$ where the function $K:=\alpha(X)$ is homogeneous of degree $1$ in $p$.
\end{proposition}
Hamiltonians $K: \T^*\Z \to \mR$ that are homogeneous of degree $1$ in $p$, and their corresponding Hamiltonian vector fields $X_K$ will be simply called {\it homogeneous}. It turns out that any homogeneous Hamiltonian vector field projects to a \emph{contact vector field} on the thermodynamic phase space $\mP(T^*\Z)$. Conversely any contact vector field on $\mP(T^*\Z)$ is the projection of a homogeneous Hamiltonian vector field on $\T^*\Z$; see \cite{entropy} for details. 

Since the state properties of the thermodynamic system are specified by a Legendre submanifold $L \subset \mP(T^*\Z)$ (Gibbs' fundamental relation) or its corresponding homogeneous Lagrangian submanifold $\cL \subset \T^*\Z$, any dynamics of a thermodynamic system should leave $\cL$ invariant. This is elegantly characterized as follows.
\begin{proposition}
A homogeneous Hamiltonian vector field $X_K$ leaves a homogeneous Lagrangian submanifold $\cL$ invariant if and only if $K$ is zero on $\cL$.
\end{proposition}
\begin{remark} A similar statement holds for the corresponding Legendre submanifold $L$: the contact vector field leaves $L$ invariant if and only if its contact Hamiltonian is zero on $L$.
\end{remark}
Finally we mention, see again \cite{entropy} for details, that the Poisson bracket of two homogeneous Hamiltonian functions is homogeneous, and that the Lie bracket of two homogeneous Hamiltonian vector fields on $\T^*\Z$ is homogeneous. This allows to set up a Lie-algebraic theory for verifying \emph{controllability} and \emph{observability} \cite{toulouse, entropy} for port-thermodynamic systems.

\section{Sidebar: Digression on Optimal Control and Homogeneous Symplectic Geometry}
Recently it was shown in \cite{ohsawa, respondek} how Pontryagin's Maximum principle naturally leads to contact geometry. This stems from the fact that the vector of co-state variables in the Mayer formulation of an optimal control problem is a {\it separating vector}, which can be arbitrarily scaled; thus leading to a contact-geometric formulation. The purpose of this digression is to highlight the initial homogeneous symplectic formulation, already implicitly present in \cite{ohsawa, respondek}. This sheds additional light on the {\it abnormal case} in optimal control.

Consider the optimal control problem of minimizing for fixed final time $\tau$ the cost criterion
\bq
\int_{0}^{\tau} L(x(t),u(t)) dt, \quad  \; x \in \mR^n,
\eq
over all input functions $u:[0,\tau] \to \mR^m$ for the dynamics $\dot{x}=f(x,u)$ with given initial condition $x(0)$.
First step is to define an additional state variable $x_0$ such that $\dot{x}_0=L(x,u), \; x_0(0)=0$. This converts the optimal control problem into the 'Mayer problem' of minimizing $x_0(T)$ over the augmented dynamics
\bq
\label{augmented}
\begin{array}{rcll}
\dot{x}_0 & = & L(x,u), \quad & x_0(0)=0 \\[2mm]
\dot{x} & = & f(x,u), \quad & x(0) \mbox{ given }
\end{array}
\eq
Next, define the pseudo-Hamiltonian $H: \T^* \mR^{n+1}\times \mR^m \to \mR$ as
\bq
H(x_0,x,\lambda_0,\lambda,u) = \lambda^\top f(x,u) + \lambda_0L(x,u) 
\eq
By construction $H$ is homogeneous of degree $1$ in $(\lambda_0,\lambda) \in \mR^{n+1}$. 
The corresponding homogeneous Hamiltonian vector field $X_H$ (parametrized by $u$) is
\bq
\label{Mayer}
\begin{array}{rcl}
\dot{x}_0 & = & L(x,u)\\[2mm]
\dot{x} & = & f(x,u)\\[2mm]
\dot{\lambda}_0 & = & 0\\[2mm]
\dot{\lambda}^\cT & = & -\lambda^\top \frac{\partial f(x,u)}{\partial x} -\lambda_0 \frac{\partial L(x,u)}{\partial x}
\end{array}
\eq
Obviously, in the first two lines of \eqref{Mayer} the original augmented dynamics \eqref{augmented} is recovered. It follows (since $H$ does not depend on $x_0$) that $\lambda_0$ is \emph{constant}, with $\lambda_0=0$ corresponding to the so-called \emph{abnormal} case. Here $(\lambda_0,\lambda) \in \mR^{n+1}$ should be understood as a vector of {\it homogeneous coordinates} for the cotangent spaces to the state space manifold with coordinates $(x_0,x)$. 
Hence for $\lambda_0\neq0$ the standard co-state variables are defined as
\bq
p := \frac{\lambda}{-\lambda_0},
\eq
resulting in the differential equations of Pontryagin's Maximum principle
\bq
\label{Pontryagin}
\begin{array}{rcl}
\dot{x}_0 & = & L(x,u)\\[2mm]
\dot{x} & = & f(x,u)\\[2mm]
\dot{p}^\cT & = & -p^\top \frac{\partial f}{\partial x}(x,u) + \frac{\partial L}{\partial x^\top}(x,u)
\end{array}
\eq
(where the equation for $\dot{x}_0$ could be left out as well). As noted in \cite{ohsawa,respondek}, the dynamics \eqref{Pontryagin} is a {\it contact vector field} on the odd-dimensional contact manifold with coordinates $x_0,x,p$ and the contact form $dx_0 - \sum_{i=1}^np_idx_i$. Alternatively, regarded from a homogeneous symplectic point of view, \eqref{Pontryagin} is the {\it projection} of the homogeneous Hamiltonian vector field \eqref{Mayer}.
Thus the differential equations of Pontryagin's Maximum principle can be understood as arising from the choice of the 'intensive' variables $p= \frac{\lambda}{-\lambda_0}$ in case the constant $\lambda_0$ is different from zero (the 'normal case'). On the other hand, in the {\it abnormal} case $\lambda_0=0$ the differential equations of the Maximum principle are most easily given in the form
\bq
\label{abnormal}
\begin{array}{rcl}
\dot{x}_0 & = & L(x,u)\\[2mm]
\dot{x} & = & f(x,u)\\[2mm]
\dot{\lambda}^\cT & = & -\lambda^\top \frac{\partial f}{\partial x}(x,u),
\end{array}
\eq
and thus do {\it not} correspond to a choice of 'intensive variables' (such as, for example, $\widetilde{p}= (\lambda_0, \lambda_2, \cdots, \lambda_n)^\top\cdot\frac{1}{-\lambda_1}$ in case $\lambda_1 \neq 0$.) This explains the peculiar form of the differential equations of the Maximum principle in the abnormal case.

Note furthermore that for the \emph{infinite-horizon} ($\tau \to \infty$) optimal control problem, the stationary Hamilton-Jacobi-Bellman equation corresponds to a homogeneous Lagrangian submanifold $\cL \subset T^* \mR^{n+1}$, with generating function $-\lambda_0 V(x)$, where $V$ is Bellman's \emph{value function}, i.e.,
\bq
\cL= \{(x_0,x,\lambda_0,\lambda) \mid x_0=V(x), \lambda=-\lambda_0\frac{\partial V}{\partial x}(x) \} ,
\eq
while the minimum of $H$ on $\cL$ with respect to $u$ is zero:
\bq
\min_u H(V(x),x, \lambda_0, -\lambda_0\frac{\partial V}{\partial x}(x),u)=0
\eq

\end{document}